\RequirePackage{etoolbox}
\csdef{input@path}{{style/}{graphics/}}
\documentclass[MSNbibl,number,citesort,seceqn,dvips]{arxbj}
\usepackage{upgreek}
\usepackage{graphicx}

\aid{0}
\volume{21}
\issue{1}
\pubyear{2015}
\firstpage{437}
\lastpage{464}
\doi{10.3150/13-BEJ574} 

\makeatletter

\newcommand{\rrVert}{\Vert}
\newcommand{\rrvert}{\vert}
\newcommand{\llVert}{\Vert}
\newcommand{\llvert}{\vert}

\newtheorem{theorem}{Theorem}

\newtheorem{coro}{Corollary}

\newremark{rema}{Remark}
\newtheorem{lema}{Lemma}

\newcommand{\E}{{\mathbf E}}
\newcommand{\ls}{\leq}
\newcommand{\gs}{\geq}

\newcommand{\TV}{\operatorname{TV}}

\def\sfrac#1#2{#1/#2}
\def\vfrac#1#2{(#1)/#2}
\def\afrac#1#2{#1/(#2)}
\def\vafrac#1#2{(#1)/(#2)}

\renewcommand{\emptyset}{\varnothing}
\makeatother

\begin{document}
\begin{frontmatter}

\title{Integrability and concentration of the truncated variation for
the sample paths of fractional Brownian motions, diffusions and L\'
{e}vy processes}

\runtitle{Integrability of the truncated variation}

\begin{aug}
\author[A]{\inits{W.M.}\fnms{Witold Marek} \snm{Bednorz}\thanksref{A}\ead[label=e1]{wbednorz@mimuw.edu.pl}}
\and
\author[B,C]{\inits{R.M.}\fnms{Rafa\L{} Marcin} \snm{\L{}ochowski}\corref{}\thanksref{B,C}\ead[label=e2]{rlocho314@gmail.com}\ead[label=u1,url]{http://akson.sgh.waw.pl/\textasciitilde rlocho}}
\address[A]{Institute of Mathematics, Warsaw University, Banacha 2,
02-097 Warszawa, Poland.\\ \printead{e1}}
\address[B]{Department of Mathematics and Mathematical Economics,
Warsaw School of Economics, Madali\'{n}skiego 6/8, 02-513 Warszawa, Poland.}
\address[C]{Department of Core Mathematics and Social Sciences, Prince
Mohammad Bin Fahd University, P.O. Box 1664, Al Khobar 31952, Saudi Arabia.\\
\printead{e2,u1}}
\end{aug}

\received{\smonth{11} \syear{2012}}
\revised{\smonth{9} \syear{2013}}

%
\begin{abstract}
For a real c\`{a}dl\`{a}g function $f$ defined on a compact interval,
its truncated variation
at the level $c>0$ is the infimum of total variations of functions
uniformly approximating $f$
with accuracy $c/2$ and (in opposite to the total variation) is always
finite. In this paper,
we discuss exponential integrability and concentration properties of
the truncated variation of
fractional Brownian motions, diffusions and L\'{e}vy processes. We
develop a special technique
based on chaining approach and using it we prove Gaussian concentration
of the truncated variation
for certain class of diffusions. Further, we give sufficient and
necessary condition for the
existence of exponential moment of order $\alpha>0$ of truncated
variation of L\'{e}vy process
in terms of its L\'{e}vy triplet.
\end{abstract}

%
\begin{keyword}
\kwd{diffusions}
\kwd{Gaussian processes}
\kwd{L\'{e}vy processes}
\kwd{sample boundedness}
\kwd{truncated variation}
\end{keyword}

\end{frontmatter}
%
\section{Introduction}

Let $X= ( X(t) )_{t \geq0}$ be a real valued
stochastic process with c\`{a}dl\`{a}g trajectories. In general, the total
path variation of $X$ on the compact interval $[a,b] \subset[0,+\infty
)$, defined as
\begin{eqnarray*}
\TV \bigl( X,[a,b] \bigr) =\sup_{n}\sup_{a\leq
t_{0}<t_{1}<\cdots<t_{n}\leq b}
\sum_{i=1}^{n}\bigl\llvert
X(t_{i})-X(t_{i-1})\bigr\rrvert ,
\end{eqnarray*}
may be (and in many most important cases is) almost surely infinite.
However, in the neighborhood of every c\`{a}dl\`{a}g path we may easily
find a function with finite total variation.

Let $f$ be a c\`{a}dl\`{a}g function $f\dvtx [a,b] \rightarrow
\mathbb{R} $ and let $ c>0$. The natural question arises, what is the
smallest possible (or the greatest lower bound for the) total variation
of functions from the
ball $B(f,c/2)= \{ g\dvtx \llVert  f-g\rrVert  _{\infty}\leq
c/2 \} $,
where $\llVert  f-g\rrVert  _{\infty}:=\sup_{s\in[a,b] }\llvert
f ( s ) -g ( s ) \rrvert  $. Some bound from
below reads as
\begin{eqnarray*}
\TV \bigl( g,[a,b] \bigr) \geq \TV^{c} \bigl( f,[a,b] \bigr) ,
\end{eqnarray*}
where
%
\begin{equation}\label{tv:def}
\TV^{c} \bigl( f,[a,b] \bigr) :=\sup_{n}\sup
_{a\leq
t_{0}<t_{1}<\cdots<t_{n}\leq b}\sum_{i=1}^{n}\max
\bigl\{ \bigl\llvert f ( t_{i} ) -f ( t_{i-1} ) \bigr
\rrvert -c,0 \bigr\}
\end{equation}
and follows immediately from the inequality
\begin{eqnarray*}
\bigl\llvert g ( t_{i} ) -g ( t_{i-1} ) \bigr\rrvert \geq
\max \bigl\{ \bigl\llvert f ( t_{i} ) -f ( t_{i-1} ) \bigr
\rrvert -c,0 \bigr\} .
\end{eqnarray*}
It is possible to show (cf. \L ochowski \cite{Loch}) that in
fact we have equality
%
\begin{equation}\label{tv2c}
\inf \bigl\{ \TV \bigl( g,[a,b] \bigr) \dvtx \llVert f-g\rrVert _{\infty}\leq
c/2 \bigr\} =\TV^{c} \bigl( f,[a,b] \bigr)
\end{equation}
attained for some function $f^c$ from the ball $B(f,c/2)$.

\begin{rema}
Since we deal with c\`{a}dl\`{a}g functions, a more natural setting of
our problem would be the investigation of
\begin{eqnarray*}
\inf \bigl\{ \TV \bigl( g,[a,b] \bigr) \dvtx  g-  \mbox{c\`{a}dl\`{a}g},
d_{D}(f,g)\leq c/2 \bigr\},
\end{eqnarray*}
where $d_{D}$ denotes the Skorohod metric. Since the total variation
does not depend on the (continuous and strictly increasing) change of
argument and the function $f^c$ minimizing $ \TV ( g,[a,b]
)$ appears to be a c\`{a}dl\`{a}g one, solutions of both problems coincide.
\end{rema}

The quantity (\ref{tv:def}) is called truncated variation and it is
finite for any c\`{a}dl\`{a}g function, since every such a function may
be uniformly approximated by step functions. Moreover, the truncated
variation is a
continuous and convex function of the parameter $c>0$ (cf. \L ochowski \cite{Loch}) and it obviously tends to the total variation as
$c \downarrow
0$. For a process with paths with almost surely infinite total
variation may be of interest to assess the rate at which $\TV^c$
diverges to infinity.

This was done so far for continuous semimartingales and it appears (cf.
\L ochowski and Mi\l o\'{s} \cite{LM}) that for any continuous semimartingale $X$
we have that
%
\begin{equation}
\label{LLN} c \cdot \TV^{c} \bigl( X,[a,b] \bigr)
\rightarrow_{c\downarrow0} \langle X \rangle_{b}- \langle X
\rangle_{a} \qquad \mbox{almost surely},
\end{equation}
where $ \langle\cdot \rangle$ denotes the quadratic variation
of $X$. The truncated variation appears also implicitly in the paper
Picard \cite{P08} where it corresponds to the double
Lebesgue measure $L^c$
of a trimmed tree at the level $c$, associated with a c\`{a}dl\`{a}g
path. In Picard \cite{P08} there were established deep
connections of this
measure, the variation index and the upper box (or Minkowski)
dimension, as well as the counterparts of (\ref{LLN}) in terms of
$L^c$ for fractional Brownian motions and stable L\'{e}vy processes.

For $t\geq0$ denote $\TV^{c} ( X,t  ) = \TV^{c} (
X, [ 0,t ]  )$. For $X$ being the unique strong
solution of the equation
$X_0 =0, \mathrm{d}X_{t}=\mu ( X_{t} )\, \mathrm{d}t+\sigma ( X_{t}
) \,\mathrm{d}W_{t},t\in[0, S]$,
driven by a standard Brownian motion $W$, with $\mu$ and $\sigma$
satisfying some linear growth conditions, we have second order
convergence result (cf. \L ochowski and Mi\l o\'{s} \cite{LM}, Theorem~10)
%
\begin{equation}
\label{CLT} \TV^{c} ( X,t ) -\frac{ \langle
X \rangle_{t}}{c} \Rightarrow_{c\downarrow0}
\tilde{W}_{
\langle X \rangle_{t}/3},
\end{equation}
where $\tilde{W}$ is a standard Brownian motion independent from $W$
and the convergence ``$\Rightarrow$'' is understood as the weak
functional convergence in $C([0, S],{\mathbb R})$ topology.

The truncated variation is more informative than $p$-variation, since
the latter may be described in terms of the asymptotic properties of
$\TV^c$ as $c \downarrow0$ but for any fixed $c>0$, $\TV^c(X,S)$ is a
proper random variable and it is possible to consider its distribution.
For $X = W$ and fixed \mbox{$S,c>0$} convergence result (\ref{CLT}) seems to
indicate very strong concentration of $\TV^c(W,S)$ around $S/c$, but it
still does not tell anything about the tail probabilities of the
functional considered.

These observations motivated us to study the integrability and
concentration properties of the truncated variation in greater detail.
Some investigation into this direction was already undertaken in
\L ochowski \cite{Loch1}, where the existence of the moment
generating function of the
truncated variation of Brownian motion with drift on the whole real
line was proven. In this paper, we obtain much stronger -- Gaussian
concentration result, by which we mean the integrability of $\exp
(\alpha \TV^c(X,S)^2  )$ for some positive $\alpha$.

Another incentive for the study of the magnitude of truncated variation
for possibly broad class of processes is the pathwise approach to
stochastic integration. In \L ochowski \cite{Loch2}, it was
shown that when both --
integrand and integrator are semimartingales then it is possible to
define the stochastic integral, with some correction term, as an almost
sure limit of pathwise Lebesgue--Stieltjes integrals. The construction
utilizes uniform approximation of the integrator with finite variation
processes. The truncated variation gives the magnitude of such
integrals, more precisely
\begin{eqnarray*}
\inf_{\llVert  X-X^c\rrVert _{\infty}\leq c/2} \sup_{\llVert  Y\rrVert _{\infty}\leq1}\int
_{0}^{S}Y_{-}\,\mathrm{d} X^{c}= \inf
_{\llVert  X-X^c\rrVert _{\infty}\leq c/2} \TV \bigl(X^c,S \bigr)= \TV^{c}
(X,S ),
\end{eqnarray*}
where the supremum is over all c\`{a}dl\`{a}g processes $Y$ with
absolute value uniformly bounded by $1$ and the infimums are over all
pathwise c\`{a}dl\`{a}g approximations $X^c$ of $X$ such that $\llVert  X-X^c\rrVert _{\infty} := \sup_{t \geq0} \llvert
X(t)-X^c(t) \rrvert \leq c/2$.

In this paper, we study the magnitude of the truncated variation for a
broad class of stochastic processes, including Gaussian processes,
among them fractional Brownian motions, and diffusions. Further we also
consider L\'{e}vy processes.
Our main goal is to describe the tail behavior of $\TV^c(X,S)$
assuming that $X$ satisfies some increment condition.
We use various techniques depending on the assumption we make.

At the beginning, we use the chaining concept, we assume that $X$
satisfies some exponential integrability condition on increments and deduce
the exponential integrability of the truncated variation (e.g.,
diffiusions with bounded covariance and drift coefficients).
The chaining approach was first used to study problems of sample
boundedness of processes on the general index space Fernique
\cite{Fer1,Fer2}.
The method was developed
to give the full description of classes of processes that are sample
bounded, under certain integrability condition
Bednorz \cite{Bed1,Bed2},
Bednorz \cite{Bed3}, Ledoux and Talagrand \cite{Le-Ta},
Talagrand \cite{Tal1}, and the small ball probability
Li and Shao \cite{Li}.
For a comprehensive study where many analytical examples are given, see
Talagrand \cite{Tal2}.
In our study, we need some modification of this idea, since we are
interested in bounding
the supremum of special sums of increments, not the supremum over
increments itself.
Therefore, we have to invent a special random variable of exponential
integrability that bounds the truncated variation.

Our main guiding example is the class of fractional Brownian motions,
that is, centered Gaussian processes $W_H$, $H\in(0,1)$, starting from
$0$ and such that $\mathbf{E} \llvert W_H(t)-W_H(s)\rrvert ^2=|t-s|^{2H}$.
One of the corollaries we get is the following concentration inequality
\begin{eqnarray*}
\mathbf{P} \bigl(\TV^c(W_H,S)\geq c^{\vfrac{H-1}{H}}S(A_H+B_Hu)
\bigr)\leq C_H\exp\bigl(-u^{2H}\bigr),\qquad  \mbox{for } u\geq0,
\end{eqnarray*}
where $A_H,B_H,C_H$ are constants; moreover, for $H\geq\frac{1}{2}$
one can set $C_H=1$.
By the homogeneity of increments, we deduce that for $Sc^{-1/H} \geq
2$, $\mathbf{E} \TV^c(W_H,S)$ is comparable with $c^{\vfrac{H-1}{H}}S$ and
in this way we prove that for $u\geq0$,
%
\begin{equation}
\label{hhh} \mathbf{P} \bigl(\TV^c(W_H,S)\geq\mathbf{E}
\TV^c(W_H,S) (\bar {A}_H+
\bar{B}_Hu) \bigr)\leq\bar{C}_H\exp\bigl(-u^{2H}
\bigr),
\end{equation}
for some constants $\bar{A}_H,\bar{B}_H,\bar{C}_H$ (again $\bar
{C}_H=1$ for $H\geq\frac{1}{2}$).
In fact, any process with similar properties as the fractional Brownian
motion, that is, satisfying some boundedness condition of the
increments (inequality (\ref{cond1})) may be treated by our method.

Next, we turn to investigate the standard Brownian motion, that is,
$W=W_{1/2}$, and diffusions driven by it. Here we can improve our
result using
the Markov property. It turns out that for Markov processes with
moderate growth some local exponential integrability can be extended to
the global one. Note that (\ref{hhh}) implies
the existence of the Laplace transform $\mathbf{E}\exp (\alpha
\TV^c(W,S)  )$ for sufficiently small $\alpha> 0$; assuming the
Markov property for diffusions with moderate growth we get the estimate
for the Laplace transform of their truncated variations on the whole
real line.
The main result we get this way is Theorem~\ref{diffusions}, which for
a standard Brownian motion and $Sc^{-2} \geq2$ implies the following
concentration inequality
\begin{eqnarray*}
\mathbf{P} \bigl(\TV^c(W,S)\geq\bar{A} \mathbf{E}
\TV^c(W,S)+ \bar {B} \sqrt{S}u \bigr)\leq\exp\bigl(-u^2
\bigr), \qquad \mbox{for } u\geq0,
\end{eqnarray*}
here $\bar{A}, \bar{B}$ are universal constants.
Therefore, the Gaussian concentration holds for the truncated variation
of the standard Brownian motion.
Our result gives better understanding of the already mentioned result
(\ref{CLT}) from which follows that $S^{-\sfrac{1}{2}}(\TV^c(W,S)-S/c)$
converges in distribution to $\mathcal{N}(0,1/3)$ as $c \downarrow0$.

We conclude the paper by proving sufficient and necessary condition for
the finiteness of $\mathbf{E} \exp ( \alpha \TV^c(X,S)  )$
for a L\'{e}vy process $X$, in terms of its generating triplet.
Here we apply the method of level crossing stopping times.

The structure of the paper is as follows. In Section~\ref{sect2}, we
introduce the chaining approach which will lead
us to the main result on the concentration for processes with
increments of exponential decay.
Then in Section~\ref{sect3}, we discuss the application of the
developed methodology to the fractional Brownian motions and then, in
Section~\ref{sect4} its improvement for a standard Wiener process and
diffusions with moderate growth. In Section~\ref{sect5} we deal with
truncated variation of L\'{e}vy processes.

\begin{rema}
In the whole paper, any dependence of a nonnegative constant on some
parameters is always indicated by listing them in brackets or in
subscripts, for example, $C(n, \varepsilon)$ or $C_{n, \varepsilon}$.
\end{rema}

\section{The chaining approach}\label{sect2}

In this section, we prove the fundamental Theorem~\ref{Fund},
which will allow us to establish integrability and concentration
properties of the truncated variation for a broad class of processes
satisfying some increment condition.

For simplicity, we consider processes indexed by a parameter from the
metric space $(T,d)$, where $T$ is the compact interval $[0,S]$, $S>0$,
equipped with the distance
$d(s,t)=|s-t|^{q}$, $s,t\in T$, where $0<q< 1$.
Further, we introduce an Orlicz function $\varphi\dvtx  [0, +\infty)
\rightarrow\mathbb{R}$ --  convex, even, satisfying $\varphi(0)=0$,
$\varphi(1)=1$, strictly increasing and such that there exists
$L<+\infty$ such that for any $x,y\geq0$,
%
\begin{equation}
\label{lpchwyt} \varphi^{-1}(xy)\leq L\bigl(\varphi^{-1}(x)+
\varphi^{-1}(y)\bigr).
\end{equation}
Moreover, we will require that $ x \mapsto\varphi (x^q )$,
$x\geq0$, is also convex.

\begin{rema}\label{felek1}
The convexity assumptions of $\varphi$ may be weakened in such a way
that $\varphi$ is convex on some interval $[C_{\varphi},\infty)$,
where $C_{\varphi} \geq0$, and $\varphi (x^q )$ is convex
on some interval $[C_{\varphi, q},\infty)$, where $C_{\varphi,q}
\geq0$.
\end{rema}

The standard example of functions with properties mentioned are
$\varphi_p(x)=2^{x^p}-1$, $p>0$, for which condition (\ref{lpchwyt})
holds with $L_p=\max\{1,2^{\vfrac{1-p}{p}}\}$. Note that when $p\gs
1$, $\varphi_p$ is convex on whole interval $[0, +\infty)$ but when $0<p<1$,
$\varphi_p$ is convex only on the interval $[C_p;+\infty)$ where
$C_p=  (\frac{1-p}{p \ln2}  )^{1/p}$.
Clearly $\varphi_p(x^q) = \varphi_{pq}(x)$ and therefore this function
is convex on the whole interval $[0; +\infty)$ if $pq\geq1$ and
convex on the interval $[C_{p,q}; +\infty)$, where $C_{p,q}=C_{pq}$,
if $pq<1$.
We use the notation $C_p, C_{p,q}$ for all $p>0$, $0<q<1$, setting
$C_p=0$ for\vspace*{2pt} $p \geq1$ (thus $C_{p,q}=0$ for $pq\geq1$). Further, we
denote $D_p = \varphi_p(C_p)$, $D_{p,q}=\varphi_p(C_{p,q}^q)=\varphi
_{pq}(C_{pq})$.
Note that $D_{p,q}=0$ for $pq \geq1$. In more general case, we will
denote $D_{\varphi} = \varphi (C_{\varphi}  )$ and
$D_{\varphi,q} = \varphi (C_{\varphi,q}^q  )$.

Let now $X(t)$, $t\in T$, be a stochastic process with increments
controlled by $\varphi$. Namely
%
\begin{equation}
\label{cond1} \mathbf{E} \varphi \biggl(\frac{|X(s)-X(t)|}{Cd(s,t)} \biggr)\leq1
\end{equation}
for $s,t\in T$, $s \neq t$, where $0< C<\infty$ is a universal constant.

\begin{rema}\label{felek2}
In fact, in (\ref{cond1}) one may consider any distance $d$ of the
form $d(s,t)=\eta(|s-t|)$, where $\eta$ is positive, concave,
increasing to $\infty$ and such that
$\eta(0)=0$. We choose $\eta(x)=x^q$, $0<q<1$, for the sake of
simplicity, however we stress that
our results can be easily extended to a more general $\eta$.
\end{rema}

Condition (\ref{cond1}) enables us to control the magnitude of the
increments of the process $X$, while the truncated variation takes into
account only increments greater than $c$ (cf. formula (\ref{tv:def})).
Note that as the consequence of (\ref{cond1}) and the compactness of
$T$ we obtain the existence of a separable modification of $X(t)$,
$t\in T$. Then by the linear order of $T$ we can define the c\`{a}dl\`
{a}g modification of $X$
which we refer to from now on.

The fundamental result of this paper, from which exponential
integrability and concentration properties will follow, is the
following theorem.

\begin{theorem}
\label{Fund}
Let $X(t)$, $t\in T$, satisfies (\ref{cond1}). Then there exist random
variables $Z_1,Z_2\geq0$ such that $\mathbf{E} Z_1,\mathbf{E}
Z_2\leq1$ and for some universal constants $K_1(q),K_2(\varphi
,q)<\infty$ the following estimate holds
\[
\TV^c(X,S)\leq c^{\vfrac{q-1}{q}}S \bigl[K_1(\varphi,q)
\varphi ^{-1}(Z_1+D_{\varphi})+K_2(
\varphi,q) \bigl[\varphi ^{-1}(Z_2+D_{\varphi,q})
\bigr]^{\sfrac{1}{q}} \bigr].
\]
\end{theorem}

\begin{rema}
The main reason why the result holds is that (\ref{cond1}) gives an
exponential decay of increments with large jumps. Therefore,
we can show a global upper bound on increments in the defined set
approximation of the truncated variation. Such an
idea is used to bound suprema of processes, for example, Bednorz \cite{Bed1}, Fernique \cite{Fer1},
Kwapie\'n and Rosi\'nski \cite{Kwap} and Talagrand \cite{Tal1}. In this paper,
the main technical contribution is to invent a common upper bound for
an arbitrary sum of truncated increments.
\end{rema}

The meaning of the result the that for suitable $\varphi$ and $0<q<1$
there holds some concentration inequality.
To formulate results in an elegant way, observe that there exists $
E_{q} \in[0;1]$ such that
$E_q+x^{1/q}\geq x$ for $x\gs0$ and hence we get
%
\begin{equation}\label{e}
E_{q}+ \bigl[\varphi^{-1} \bigl(x+\max \{D_{\varphi},
D_{\varphi,q} \} \bigr) \bigr]^{\sfrac{1}{q}}\geq\varphi ^{-1}(x+D_{\varphi})
\qquad \mbox{for } x\geq0.
\end{equation}
As a consequence of Theorem~\ref{Fund}, (\ref{e}) and Jensen's
inequality we get the following corollary.

\begin{coro}\label{cord1}
Under the assumptions of Theorem~\ref{Fund} there exist r.v. $Z$ such
that $Z\geq0$, $\mathbf{E} Z\leq1$ and for some constants
$A_{\varphi,q}$, $B_{\varphi,q}$ the following estimate holds
\[
\TV^c(X,S)\leq c^{\vfrac{q-1}{q}}S \bigl[A_{\varphi,q}+B_{\varphi
,q}
\bigl[\varphi^{-1}\bigl(Z+\max \{D_{\varphi}, D_{\varphi
,q}
\}\bigr) \bigr]^{\sfrac{1}{q}} \bigr].
\]
\end{coro}

For $\varphi= \varphi_p$ let us denote $A_{p,q} = A_{\varphi,q}$ and
$B_{p,q} = B_{\varphi,q}$. Applying Corollary~\ref{cord1}, the Markov
inequality and the fact that $D_{p,q}\geq D_{p}$ we obtain:

\begin{coro}\label{cord2}
Let $X(t)$, $t\in T$, satisfies (\ref{cond1}) with $\varphi= \varphi
_p$. The following inequality holds
\[
\mathbf{P} \bigl(\TV^c(X,S)\geq c^{\vfrac{q-1}{q}}S[
\bar{A}_{p,q}+\bar {B}_{p,q}u] \bigr)\leq\bar{D}_{p,q}
\exp\bigl(-u^{pq}\bigr), \qquad \mbox{for } u>0,
\]
where $\bar{A}_{p,q},\bar{B}_{p,q}$ are universal constants, $\bar
{A}_{p,q}=A_{p,q}+(2/\ln2)^{\afrac{1}{pq}}B_{p,q}$, $\bar
{B}_{p,q}=\linebreak[4] (2/\ln2)^{\afrac{1}{pq}}B_{p,q}$ and $\bar{D}_{p,q}=D_{p,q}+1$.
In particular, $\bar{D}_{p,q}=1$ for $pq\geq1$.
\end{coro}

To prove Theorem~\ref{Fund}, we start with the construction of finite
sets approximating $T$.

\subsection{Approximating sequence}

The first tool we need is a proper geometric approximation of the set
$T$. The approximation consists of a sequence of finite sets
$(T_n)^{\infty}_{n=0}, T_n \subset T$
constructed in such a way that for each point $t\in T$ and $n=0,
1,2,\ldots\,$, there exists a point $s\in T_n$, such that $s\ls t$ and
$d(s,t)\leq r^{-n}S^q$. Here, we fix $r\gs4$.
One of possible constructions is the following
%
\begin{equation}
\label{T_n} T_n= \bigl\{kr^{-\sfrac{n}{q}}S\dvtx  k =0, 1, 2, \ldots
\bigr\}\cap T.
\end{equation}
For $T_n$ defined by (\ref{T_n}) and $t \in T$, by $\pi_n(t)$ we
denote the unique point $s\in T_n$ such that $s\ls t$ and $d(s,t)<
r^{-n}S^q$. This way we define the function $\pi_n\dvtx T\rightarrow T_n$.
We have $d(t,\pi_n(t)) < r^{-n}S^q$ for all $t\in T$
and $\pi_n(s)\ls\pi_n(t)$ if $s\ls t$. Note also that for $s,t\in
T_n$, $s\neq t$, $d(s,t)\gs r^{-n}S^q$.
Clearly
%
\begin{equation}
\label{ham0} r^{\sfrac{n}{q}} <|T_n|= \bigl\lfloor r^{\sfrac{n}{q}}
\bigr\rfloor+1 \ls r^{\sfrac{n}{q}} +1.
\end{equation}
Moreover for any $m=1,2,\dots$
%
\begin{equation}
\label{ham1} \sum^m_{n=0}r^{-n}|T_{n+1}|
\leq\sum^m_{n=0}r^{-n}
\bigl(r^{\vfrac{n+1}{q}}+1\bigr)\leq A(r,q) r^{m\vfrac{1-q}{q}},
\end{equation}
where $A(r,q):=r^{\vfrac{2-q}{q}}(r^{\vfrac{1-q}{q}}-1)^{-1}$ (note
that $r \geq2$).
For each $t\in T_{n+1}$ let $I_{n+1}(t)$ denote
the set of the nearest neighbors of $t$ in $T_{n+1}$, namely
%
\begin{equation}
I_{n+1}(t) = \bigl\{s\in T_{n+1} \dvtx  d(s,t)
\leq2r^{-n}S^q \bigr\}.
\end{equation}
Observe that
since $|s-t|\geq r^{-\vfrac{n+1}{q}}S$ for $s,t\in T_{n+1}$, $s\neq t$,
%
\begin{equation}
\label{ham2} \bigl|I_{n+1}(t)\bigr|\leq\frac{2^{\sfrac{1}{q}}r^{-\sfrac{n}{q}}S}{r^{-\vfrac
{n+1}{q}}S}+1=2^{\sfrac{1}{q}}r^{\sfrac{1}{q}}+1=:B(r,q).
\end{equation}

\subsection{Proof of the main theorem}

The plan of the proof is the following. After having constructed the
set approximation of $T$, we use this approximation to build a type of
discretization of any given partition and
derive a chaining bound on the truncated variation (Lemma~\ref{pierwszy}).
Then we turn to estimate each increment in the partition bound (Lemma~\ref{nowy2})
and finally apply the bounds as well as some technical observations
(Lemmas \ref{convexf}, \ref{convexf1} and \ref{nowy1})
to derive the required bounds (Lemmas \ref{drugi}, \ref{trzeci}).

Our first step is to analyze a given partition $\Pi_n=\{t_0,t_1,\dots
,t_n\}$, where $0\leq t_0<t_1<\cdots<t_n \leq S$.
We decompose the set $\{1,\dots,n\}$ into subsets $J_m$,
$m=0,1,2,\ldots\,$, defined in the following way
\[
J_m=\bigl\{i\in\{1,\dots,n\}\dvtx  r^{-m-1}S^q<d(t_{i-1},t_i)
\leq r^{-m}S^q\bigr\}.
\]
Let $M_0:=12CL$, where $L$ and $C$ are constants appearing in (\ref
{lpchwyt}) and (\ref{cond1}). The level $m_0 \in \{0,1,2,\ldots
\}$ such that
\[
r^{-m_0-1}S^q<c/M_0\ls r^{-m_0}S^q
\]
will be of particular meaning in the proof.
Since $\Pi_n$ is finite, $J_m=\emptyset$ for $m$ large enough, say
$m\gs N_0\gs m_0$.
We will use different bounds for $i\in J_m$ with $m> m_0$ and for $i\in
J_m$ with $m\ls m_0$. Therefore, let us make the trivial separation\vspace*{-1pt}
%
\begin{eqnarray}\label{czuma1}
\sum^n_{i=1} \bigl(\bigl|X(t_i)-X(t_{i-1})\bigr|-c
\bigr)_{+} & \leq& \sum^{m_0}_{m=0}
\sum_{i\in J_m} \bigl(\bigl|X(t_i)-X(t_{i-1}
)\bigr| - c\bigr)_+
\nonumber
\\[-8.5pt]\\[-8.5pt]
&& {}+\sum^{\infty}_{m=m_0+1}\sum
_{i\in J_m} \bigl(\bigl|X(t_i)-X(t_{i-1})\bigr|-c
\bigr)_{+}.\nonumber
\end{eqnarray}
Now we turn to describe the chaining method which is the main tool in
the proof. First, we fix $N\gs N_0$ and define $t^{N+1}_i=\pi
_{N+1}(t_i)$, then for $l\in\{0,1,\dots, N\}$ we put by the reverse
induction $t^l_i=\pi_{l}(t^{l+1}_i)$.
Note that by the construction of $\pi_l$ we preserve the order of the
projections, namely $t^l_0\leq t^l_1\leq\cdots\leq t^{l}_n$ for any
$0\leq l\leq N+1$.
Moreover since $N\gs N_0$ points $\{t^{N+1}_0,t^{N+1}_1,\dots
,t^{N+1}_n\}$ are separated, that is, $t^{N+1}_i\neq t^{N+1}_{i-1}$,
$i\in\{1,\dots,n\}$.
Let us denote $\bar{m}=\max \{m,m_0 \}$. For $i\in J_m$
with $m>m_0$, we estimate\vspace*{-1pt}
%
\begin{eqnarray}\label{czuma3}
&& \bigl(\bigl|X(t_i)-X(t_{i-1})\bigr|-c \bigr)_{+}
\nonumber
\\[-0.5pt]
&&\quad \leq \biggl(\bigl|X\bigl(t^{m+1}_{i}\bigr)-X
\bigl(t^{m+1}_{i-1}\bigr)\bigr|-\frac{c}{3}
\biggr)_{+}+\sum_{s\in\{i-1,i\}}\bigl|X
\bigl(t^{N+1}_s\bigr)-X(t_s)\bigr|
\\[-0.5pt]
&&\qquad {} +\sum^{N}_{l=m+1}\sum
_{s\in\{i-1,i\}} \biggl(\bigl|X\bigl(t^l_{s}
\bigr)-X\bigl(t^{l+1}_s\bigr)\bigr|-2^{-l+\bar{m}}
\frac{c}{3} \biggr)_{+}\nonumber
\end{eqnarray}
and for $i\in J_m$ with $m\leq m_0$ we have\vspace*{-1pt}
%
\begin{eqnarray}\label{czuma4}
&& \bigl(\bigl|X(t_i)-X(t_{i-1})\bigr|-c \bigr)_{+}\nonumber\\
&&\quad \leq \bigl|X\bigl(t^{m+1}_i\bigr)-X\bigl(t^{m+1}_{i-1}
\bigr)\bigr|
\nonumber
\\[-8pt]\\[-8pt]
&&\qquad {} +\sum_{s\in\{i-1,i\}}\bigl|X\bigl(t^{N+1}_s
\bigr)-X(t_s)\bigr| +\sum^{m_0}_{l=m+1}
\sum_{s\in\{i-1,i\}}\bigl|X\bigl(t^l_s
\bigr)-X\bigl(t^{l+1}_s\bigr)\bigr|\nonumber
\\[-0.5pt]
&&\qquad {} +\sum^{N}_{l=m_0+1}\sum
_{s\in\{i-1,i\}} \biggl(\bigl|X\bigl(t^l_{s}
\bigr)-X\bigl(t^{l+1}_s\bigr)\bigr|-2^{-l+\bar{m}}
\frac{c}{3} \biggr)_{+}.\nonumber
\end{eqnarray}
Putting together estimates (\ref{czuma1}), (\ref{czuma3}) and (\ref
{czuma4}), we obtain the following decomposition lemma.

\begin{lema}\label{pierwszy}
For any partition $\Pi_n=\{t_0,\ldots,t_n\}$, where $n\geq0$, $0\leq
t_0<t_1<\cdots<t_n\leq S$ and $N>m_0$
the following estimate holds\vspace*{-1pt}
\begin{eqnarray*}
\sum^n_{i=1}\bigl(\bigl|X(t_i)-X(t_{i-1})\bigr|-c
\bigr)_{+}  \leq V_1+V_2+W_1+W_2
+\sum^n_{i=1}\sum
_{s\in\{i-1,i\}}\bigl|X(t_s)-X\bigl(t^{N+1}_s
\bigr)\bigr|,
\end{eqnarray*}
where
\begin{eqnarray*}
V_1&:=&\sum^{m_0}_{m=0}\sum
_{i\in J_m}\sum^{m_0}_{l=m+1}
\sum_{s\in\{
i-1,i\}}\bigl|X\bigl(t^l_s
\bigr)-X\bigl(t^{l+1}_{s}\bigr)\bigr|;
\\
W_1&:=&\sum^{m_0}_{m=0}\sum
_{i\in J_m}\bigl|X\bigl(t^{m+1}_i
\bigr)-X\bigl(t^{m+1}_{i-1}\bigr)\bigr|;
\\
 V_{2}&:=&\sum^{\infty}_{m=0}\sum
_{i\in J_m}\sum^N_{l=\bar
{m}+1}
\sum_{s\in\{i-1,i\}} \biggl(\bigl|X\bigl(t^l_s
\bigr)-X\bigl(t^{l+1}_{s}\bigr)\bigr|-2^{-l+\bar
{m}}
\frac{c}{3} \biggr)_{+};
\\
W_2&:=&\sum^{\infty}_{m=m_0+1}\sum
_{i\in J_m} \biggl(\bigl|X\bigl(t^{m+1}_i
\bigr)-X\bigl(t^{m+1}_{i-1}\bigr)\bigr|-\frac{c}{3}
\biggr)_{+}.
\end{eqnarray*}
\end{lema}

\begin{figure}

\includegraphics{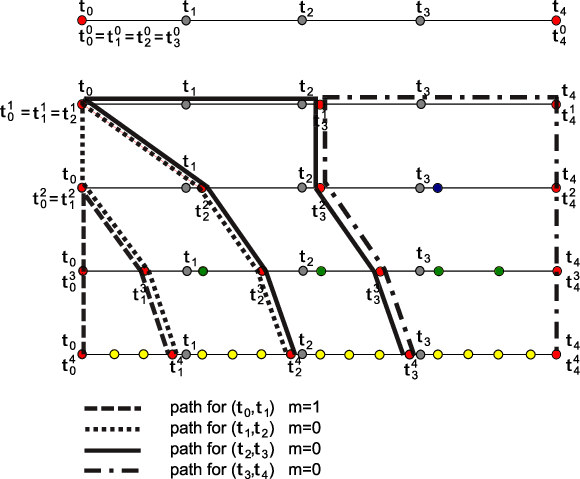}

\caption{Path approximations.}\label{fig:Path approximation}
\vspace*{-6pt}
\end{figure}

For each $i\in J_m$, $m\gs0$ we say that $t^{m+1}_s,t^{m+2}_s,\dots
,t^{N+1}_s$, $s=i-1,i$
are path approximations of $t_{i-1}$ and $t_i$, respectively (see
Figure~\ref{fig:Path approximation}). Note that for $i\in\{1,\dots,n-1\}$ there are two path
approximations of $t_i$,
one from the pair $t_{i-1},t_i$ and the second from the pair
$t_i,t_{i+1}$, that coincide, starting from some point, yet may differ
on the
length since $i\in J_m$, $i+1\in J_{m'}$ and numbers $m$ and $m'$ may
be different.
The fundamental property of the path approximation is that for a given
$u\in T_{l+1}$ the step $\pi_l(u),u$ may occur in at most two path
approximations of some $t_i$, $i\in\{0,1,\dots,n\}$.
%

\begin{lema}\label{nowy1}
Consider $u\in T_{l+1}$, $l\in\{0,1,\dots,n\}$. The step $\pi
_l(u),u$ may occur in at most two path approximations of some $t_i$,
$i\in\{0,1,\dots,n\}$, that is,
there exits no more than one $i\in\{0,1,\dots,n\}$ such that $i\in
J_m$, $m+1\ls l$ and $t^l_i=\pi_l(u)$, $t^{l+1}_i=u$ or $i+1\in
J_{m'}$, $m'+1\ls l$
and $t^l_i=\pi_l(u)$, $t^{l+1}_{i}=u$ for some $m,m'=0,1,2,\ldots,N$.
\end{lema}

\begin{pf}
Recall that $r\gs4$. It suffices to prove that for a given $i\in J_m$,
$l\gs m+1$ points $t^{l+1}_i$ and $t^{l+1}_{i-1}$
are different. Indeed since $t^{l+1}_0\ls t^{l+1}_1\ls\cdots\ls
t^{l+1}_n$ the property implies that there can be at most one $i\in\{
0,1,\dots,n\}$
such that $t^{l+1}_i=u$. To prove the assertion, we use
$d(t_i,t_{i-1})>r^{-m-1}S^q$ which implies that for $l \geq m+1$
\begin{eqnarray*}
d\bigl(t^{l+1}_i,t^{l+1}_{i-1}\bigr)
&\geq& r^{-m-1}S^q-d\bigl(t^{l+1}_{i-1},t_{i-1}
\bigr)-d\bigl(t^{l+1}_i,t_i\bigr)
\\
&\geq& r^{-m-1}S^q-2\sum^{\infty}_{j=l+1}r^{-j}S^q
\geq r^{-m-1}S^q-2\frac{r^{-m-2}S^q}{1-r^{-1}}>0.
\end{eqnarray*}
\upqed
\end{pf}
In the sequel, we will use two simple observations concerning
increasing function $\psi$ that is convex starting from some $C_0\geq
0$, that is, convex for $x\geq C_0$.

\begin{lema}
\label{convexf}
Let $\psi\dvtx  [0;+\infty) \rightarrow[0;+\infty) $ be a strictly
increasing function. Assume that $\psi$ is convex on the interval
$[C_0;+\infty)$ where $C_0\geq0$, then for any nonnegative
$x_1,\ldots,x_k$ and
positive $\alpha_1,\ldots,\alpha_k$ such that $\sum^{k}_{i=1}\alpha
_i\leq M$
we have\vspace*{-1pt}
%
\begin{equation}\label{jensen}
\sum^k_{i=1}\alpha_i
x_i \leq M\psi^{-1} \Biggl(M^{-1}\sum
^k_{i=1}\alpha_i\psi(x_i)+
\psi(C_0) \Biggr).
\end{equation}
\end{lema}

\begin{pf}
Observe that the
function $\bar{\psi}(x)=\psi (x+C_0 )-\psi(C_0)$ for
$x\geq0$ is convex, strictly increasing and such that $\bar{\psi}(0)=0$.
Consequently, $\bar{\psi}^{-1}(y)=\psi^{-1} (y+\psi(C_0)
)-C_0$ is concave with $\bar{\psi}^{-1}(0)=0$ and we have\vspace*{-1.5pt}
\begin{eqnarray*}
\sum^k_{i=1}\alpha_i
x_i &\leq&\sum^k_{i=1}
\alpha_i \psi ^{-1} \bigl(\psi(x_i)+
\psi(C_0) \bigr)
\nonumber
\\[-1pt]
& =& \sum^k_{i=1}\alpha_i
\bigl(\bar{\psi}^{-1} \bigl(\psi (x_i) \bigr)+
C_0 \bigr)\leq M C_0 + M \sum
^k_{i=1} \frac{\alpha
_i}{M} \bar{
\psi}^{-1}\bigl(\psi(x_i)\bigr)
\nonumber
\\[-1pt]
& \leq& M C_0 + M \bar{\psi}^{-1} \Biggl( \sum
^k_{i=1} \frac
{\alpha_i}{M} \psi(x_i)
\Biggr)
\nonumber
\\[-1pt]
&=&M\psi^{-1} \Biggl(M^{-1}\sum
^k_{i=1}\alpha_i\psi(x_i)+
\psi (C_0) \Biggr),\vspace*{-1pt}
\nonumber
\end{eqnarray*}
where the last inequality follows from Jensen's inequality\vspace*{-1pt}
\[
\Biggl(1 - \sum^k_{i=1}
\frac{\alpha_i}{M} \Biggr)\bar{\psi }^{-1}(0) + \sum
^k_{i=1} \frac{\alpha_i}{M} \bar{
\psi}^{-1}\bigl(\psi (x_i)\bigr) \leq\bar{
\psi}^{-1} \Biggl( \sum^k_{i=1}
\frac{\alpha
_i}{M} \psi(x_i) \Biggr).
\]
\upqed
\end{pf}
Further, we also have the following lemma.

\begin{lema}
\label{convexf1}
For any strictly increasing function $\psi\dvtx  [0;+\infty) \rightarrow
[0;+\infty) $ such that $\psi$ is convex on the interval
$[C_0;+\infty)$ where $C_0\geq0$ and for any $M >0$ and $y\geq0$, we have\vspace*{-1pt}
%
\begin{equation}\label{jensen1}
\psi^{-1} \bigl(y+\psi (C_{0} ) \bigr)\leq\max\{M,1\}
\psi^{-1} \bigl(y/M+\psi (C_{0} ) \bigr).\vspace*{-1pt}
\end{equation}
\end{lema}

\begin{pf}
Again, we consider the function $\bar{\psi}^{-1}$. If $M < 1$, then
(\ref{jensen1}) follows from the monotonicity of $\bar{\psi}^{-1}$.
Now assume that $M\geq1$.
By concavity and $\bar{\psi}^{-1} (0 )=0$, for $y\geq0$
and $M\geq1$, we get\vspace*{-1pt}
\[
M\bar{\psi}^{-1} (y/M )\geq\bar{\psi}^{-1} (y ),\vspace*{-1pt}
\]
which reads as\vspace*{-1pt}
\begin{eqnarray*}
M \bigl(\psi^{-1} \bigl(y/M+\psi (C_{0} )
\bigr)-C_{0} \bigr) & \geq& \psi^{-1} \bigl(y+\psi
(C_{0} ) \bigr)-C_{0},
\\[-1pt]
M\psi^{-1} \bigl(y/M+\psi (C_{0} ) \bigr) & \geq& \psi
^{-1} \bigl(y+\psi (C_{0} ) \bigr)+ (M-1 )C_{0}\vspace*{-1pt}
\end{eqnarray*}
and which gives\vspace*{-1pt}
\[\hspace*{100pt}
\psi^{-1} \bigl(y+\psi (C_{0} ) \bigr)\leq M\psi
^{-1} \bigl(y/M+\psi (C_{0} ) \bigr).\hspace*{100pt}\qed
\]
\noqed
\end{pf}

Now we formulate some basic bounds on increments in the chaining argument.
For simplicity, we use the following notation
\[
\Delta(u,v)=\varphi \biggl(\frac{|X(u)-X(v)|}{Cd(u,v)} \biggr),\qquad  \mbox {for all } u,v\in T.
\]
Recall that $\bar{m}=\max \{m,m_0 \}$.

\begin{lema}\label{nowy2}
Suppose that $i\in J_m$, $m\gs0$ then
\begin{enumerate}[3.]
\item[1.] for any $m\ls m_0$, $l\in\{m+1,\ldots,m_0\}$
\[
\bigl|X\bigl(t^l_i\bigr)-X\bigl(t^{l+1}_i
\bigr)\bigr|\ls Cr^{-l}S^q\varphi^{-1}\bigl(\Delta
\bigl(t^l_i,t^{l+1}_i\bigr)
\bigr);
\]
\item[2.] for any $m\gs0$, $l\in\{\bar{m}+1,\dots,N\}$
\[
\biggl(\bigl|X\bigl(t^l_i\bigr)-X\bigl(t^{l+1}_i
\bigr)\bigr|-2^{-l+\bar{m}}\frac{c}{3} \biggr)_{+}\ls
b_{l,\bar{m}} \bigl[\varphi^{-1} \bigl(a_{l,\bar
{m}}^{-1}
\Delta\bigl(t^l_i,t^{l+1}_i\bigr)
\bigr) \bigr]^{\sfrac{1}{q}},
\]
where
\[
a_{l,\bar{m}}=\varphi \biggl(\frac{2^{-l+\bar{m}}r^lc}{6CLS^q} \biggr),\qquad  b_{l,\bar{m}}=
\biggl(\frac{c}{6}2^{-l+\bar{m}} \biggr)^{\vfrac
{q-1}{q}}
\bigl(CLr^{-l}S^q \bigr)^{\sfrac{1}{q}};
\]
\item[3.] for $m\ls m_0$
\[
\bigl|X\bigl(t^{m+1}_{i}\bigr)-X\bigl(t^{m+1}_{i-1}
\bigr)\bigr|\ls2Cr^{-m}S^q\varphi^{-1}\bigl(\Delta
\bigl(t^{m+1}_i,t^{m+1}_{i-1}\bigr)
\bigr);
\]
\item[4.] and for $m>m_0$
\[
\biggl(\bigl|X\bigl(t^{m+1}_{i}\bigr)-X\bigl(t^{m+1}_{i-1}
\bigr)\bigr|-\frac{c}{3} \biggr)_{+}\ls b_m \bigl[
\varphi^{-1} \bigl(a_m^{-1}\Delta
\bigl(t^{m+1}_i,t^{m+1}_{i-1}\bigr)
\bigr) \bigr]^{\sfrac{1}{q}},
\]
where
\[
a_m=\varphi \biggl(\frac{r^mc}{12CLS^q} \biggr),\qquad  b_m=
\biggl(\frac
{c}{6} \biggr)^{\vfrac{q-1}{q}} \bigl(2CLr^{-m}S^q
\bigr)^{\sfrac{1}{q}}.
\]
\end{enumerate}
\end{lema}

\begin{rema}\label{poziomka}
Note that the choice of $M_0$ in the definition of $m_0$ guarantees that
$a_{l,\bar{m}}\gs1$ and $a_m\gs1$ for $m>m_0$. Moreover, for
$\varphi$ convex on the whole real line, that is, $C_{\varphi}=0$
one can deduce $a_{l,\bar{m}}^{-1}\ls(r/2)^{-l+\bar{m}}$ and
$a_m^{-1}\ls r^{-m+m_0+1}$.
\end{rema}

\begin{pf*}{Proof of Lemma~\ref{nowy2}}
Let us denote $u=t^l_i$, then the construction of approximation paths
implies that $t^{l+1}_i=\pi_l(u)$. Clearly $d(u,\pi_l(u))\ls
r^{-l}S^q$ and hence
%
\begin{equation}
\label{bakcyl0} \bigl|X\bigl(\pi_l(u)\bigr)-X(u)\bigr|\ls Cr^{-l}S^q
\varphi^{-1} \bigl(\Delta\bigl(\pi _l(u),u\bigr) \bigr).
\end{equation}
To prove the second assertion, we use (\ref{bakcyl0}) to get
%
\begin{equation}
\label{bakcyl1} \biggl(\bigl|X\bigl(\pi_l(u)\bigr)-X(u)\bigr|-2^{-l+\bar{m}}
\frac{c}{3} \biggr)_{+}\ls \biggl(Cr^{-l}S^q
\varphi^{-1} \bigl(\Delta\bigl(\pi_l(u),u\bigr)
\bigr)-2^{-l+m}\frac
{c}{3} \biggr)_{+}.
\end{equation}
Now we rewrite (\ref{bakcyl1}) using $a_{l,\bar{m}}$
\begin{eqnarray*}
&& \biggl(\bigl|X\bigl(\pi_l(u)\bigr)-X(u)\bigr|-2^{-l+\bar{m}}
\frac{c}{3} \biggr)_{+}
\\
&&\quad \ls \biggl[Cr^{-l}S^q \bigl(\varphi^{-1} \bigl(
\Delta\bigl(\pi _l(u),u\bigr) \bigr)-L\varphi^{-1}(a_{l,\bar{m}})
\bigr)_{+}- 2^{-l+\bar{m}}\frac{c}{6} \biggr]_{+}.
\end{eqnarray*}
Therefore, we can apply (\ref{lpchwyt}) and see
%
\begin{eqnarray}\label{bakcyl2}
&& \biggl(\bigl|X\bigl(\pi_l(u)\bigr)-X(u)\bigr|-2^{-l+\bar{m}}
\frac{c}{3} \biggr)_{+}
\nonumber
\\
&&\quad \ls \biggl(CLr^{-l}S^q
\varphi^{-1} \bigl(a_{l,\bar
{m}}^{-1}\Delta\bigl(
\pi_l(u),u\bigr) \bigr)-2^{-l+\bar{m}}\frac{c}{6}
\biggr)_{+}.
\end{eqnarray}
Using the inequality $(x-1)_{+}\ls x^{\sfrac{1}{q}}$ valid for $x\gs
0$, we get
\begin{eqnarray*}
&& \biggl(\bigl|X\bigl(\pi_l(u)\bigr)-X(u)\bigr|-2^{-l+\bar{m}}
\frac{c}{3} \biggr)_{+}
\\
&&\quad \ls\frac{c}{6}2^{-l+\bar{m}} \bigl( 6c^{-1}2^{l-\bar{m}}CLr^{-l}S^q
\varphi^{-1} \bigl(a_{l,\bar{m}}^{-1}\Delta\bigl(
\pi_l(u),u\bigr) \bigr)-1 \bigr)_{+}
\\
&&\quad \ls\frac{c}{6}2^{-l+\bar{m}} \bigl(6c^{-1}2^{l-\bar
{m}}CLr^{-l}S^q
\bigr)^{\sfrac{1}{q}} \bigl[\varphi^{-1} \bigl(a_{l,\bar{m}}^{-1}
\Delta\bigl(\pi _l(u),u\bigr) \bigr) \bigr]^{\sfrac{1}{q}}
\\
&&\quad  = b_{l,\bar{m}} \bigl[\varphi^{-1} \bigl(a_{l,\bar{m}}^{-1}
\Delta \bigl(\pi_l(u),u\bigr) \bigr) \bigr]^{\sfrac{1}{q}}.
\end{eqnarray*}
To prove the third assertion, we first observe that
$d(t_i,t_{i-1})\leq r^{-m}S^q$ for $i\in J_m$ and hence
\begin{eqnarray*}
 d\bigl(t^{m+1}_i,t^{m+1}_{i-1}\bigr)
&\leq& d(t_i,t_{i-1})+d\bigl(t^{N+1}_i,t_i
\bigr)+d\bigl(t^{N+1}_{i-1},t_{i-1}\bigr)
 +\sum^{N}_{l=m+1}\bigl[d
\bigl(t^{l+1}_i,t^{l}_i\bigr)+d
\bigl(t^{l+1}_{i-1},t^{l}_{i-1}\bigr)
\bigr]\\
&\leq& r^{-m}S^q+2\sum^{\infty}_{l=m+1}r^{-l}S^q,
\end{eqnarray*}
so
%
\begin{equation}
\label{bakcyl3} d\bigl(t^{m+1}_i,t^{m+1}_{i-1}
\bigr)\leq2r^{-m}S^q.
\end{equation}
Denoting $u=t^{m+1}_i$ and $v=t^{m+1}_{i-1}$ we get in the same way as
(\ref{bakcyl0}) that
\[
\bigl|X(u)-X(v)\bigr|\ls2Cr^{-m}S^q\varphi^{-1} \bigl(
\Delta(u,v) \bigr).
\]
Then using the same idea as for the second assertion we deduce the
remaining inequality.
\end{pf*}
We turn to apply the above lemmas to bound increments in the chaining
bound formulated in Lemma~\ref{pierwszy}. First, we consider a bound
on $V_1+W_1$.

\begin{lema}\label{drugi}
There exists a universal constant $K_1(\varphi, r,q)<\infty$ and a
random variable $Z_1\geq0$ independent from the partition $\Pi_n$,
such that $\mathbf{E} Z_1\leq1$ and for $V_1$ and $W_1$ defined in
Lemma~\ref{pierwszy} one has
\begin{eqnarray*}
& V_1+W_1\leq K_1(
\varphi,r,q)c^{\vfrac{q-1}{q}}S\varphi ^{-1}(Z_1+D_{\varphi}).
\end{eqnarray*}
\end{lema}

\begin{pf}
By Lemma~\ref{nowy1} and the first bound in Lemma~\ref{nowy2}, we get
%
\begin{eqnarray}\label{puma1}
 V_1&\leq&2\sum^{m_0}_{l=0}
\sum_{u\in T_{l+1}}\bigl|X(u)-X\bigl(\pi_l(u)\bigr)\bigr|
\nonumber
\\[-8pt]\\[-8pt]
&\leq&\mathcal{V}_1:= 2C\sum
^{m_0}_{l=0}r^{-l}S^q\sum
_{u\in T_{l+1}}\varphi^{-1}\bigl(\Delta \bigl(
\pi_l(u),u\bigr)\bigr).\nonumber
\end{eqnarray}
To bound $W_1$ we use (\ref{bakcyl3}), that is,
that $d(t^{m+1}_i,t^{m+1}_{i-1})\leq2r^{-m}S^q$ for $i\in J_m$.
Using the already defined sets $I_m(u)=\{v\in T_{m+1}\dvtx  d(u,v)\leq
2r^{-m}S^q\}$,
and the third bound in Lemma~\ref{nowy2}
%
\begin{eqnarray}\label{puma2}
 W_1&\leq&\sum^{m_0}_{l=0}
\sum_{u\in T_{l+1}}\sum_{v\in
I_{l+1}(u)}\bigl|X(u)-X(v)\bigr|
\nonumber
\\[-8pt]\\[-8pt]
&\leq&\mathcal{W}_1:=C\sum
^{m_0}_{l=0}2r^{-l}S^q\sum
_{u\in T_{l+1}}\sum_{v\in I_{l+1}}
\varphi^{-1} \bigl(\Delta(u,v) \bigr).\nonumber
\end{eqnarray}
We calculate the sum of all weights appearing in (\ref{puma1}) and
(\ref{puma2}).
By (\ref{ham2}) for each $u\in T_{m+1}$ we have $|I_{l+1}(u)|\leq
B(r,q)$ and hence, using also (\ref{ham1})
\begin{eqnarray*}
 M_1&:=&\sum^{m_0}_{l=0}r^{-l}S^q
\biggl[|T_{l+1}|+\sum_{u\in
T_{l+1}}\bigl|I_{l+1}(u)\bigr|
\biggr]
\\
&\leq& \bigl[1+B(r,q) \bigr]S^q\sum^{m_0}_{l=0}r^{-l}|T_{l+1}|\\
&\leq& A(r,q)\bigl[1+B(r,q)\bigr]r^{m_0\vfrac{1-q}{q}}S^q.
\end{eqnarray*}
Therefore by $c\leq M_0r^{-m_0}S^q$ we get $M_0r^{m_0\vfrac
{1-q}{q}}S^q\ls
M_0^{\sfrac{1}{q}}c^{\vfrac{q-1}{q}}S$ and hence
\[
M_1\leq M_0^{\sfrac{1}{q}}A(r,q)\bigl[1+B(r,q)\bigr]
c^{\vfrac{q-1}{q}}S.
\]
Using Lemma~\ref{convexf} for $\varphi$ which is convex above
$C_{\varphi}$ we get
%
\begin{equation}
\label{puma3} \mathcal{V}_1+\mathcal{W}_1\leq2 C
M_1\varphi^{-1}\bigl(Z_1+ \varphi
(C_{\varphi})\bigr)\leq K_1(r,q)c^{\vfrac{q-1}{q}}S\varphi
^{-1}\bigl(Z_1+\varphi(C_{\varphi})\bigr),
\end{equation}
where $K_1(\varphi,r,q):=2 C M_0^{\sfrac{1}{q}} A(r,q)[1+B(r,q)] $
(the dependence on $\varphi$ is through $C$ and $L$) and
\begin{eqnarray*}
 Z_1=M_1^{-1}\sum
^{m_0}_{l=0}r^{-l}S^q\sum
_{u\in T_{l+1}} \biggl(\Delta\bigl(\pi_l(u),u
\bigr)+\sum_{v\in I_{l+1}(u)}\Delta(u,v) \biggr).
\end{eqnarray*}
Obviously $Z_1\geq0$ and $\mathbf{E} Z_1\leq1$ by (\ref{cond1}) and
the definition of $M_1$.
Combining (\ref{puma1}), (\ref{puma2}) and (\ref{puma3}) we get the result.
\end{pf}
Our second goal is to prove a bound for $V_2+W_2$ in Lemma~\ref
{pierwszy} above the level $m_0$.

\begin{lema}\label{trzeci}
There exists a universal constant $K_2(\varphi,r,q)<\infty$ and a
random variable $Z_2\geq0$ independent from the partition $\Pi_n$
such that $\mathbf{E} Z_2\leq1$ and for $V_2$ and $W_2$ defined in
Lemma~\ref{pierwszy} the following inequality holds
\[
V_2+W_2\leq K_2(\varphi,r,q) \bigl[
\varphi^{-1}(Z_2+D_{\varphi
,q}) \bigr]^{\sfrac{1}{q}}.
\]
\end{lema}

\begin{pf}
First, we prove a bound for $V_{2}$. We analyze the increment
\[
\biggl(\bigl|X\bigl(t^{l+1}_i\bigr)-X\bigl(t^l_i
\bigr)\bigr|-2^{\bar{m}-l}\frac{c}{3} \biggr)_{+},\qquad  l>\bar{m}, i\in
J_m, m\gs0.
\]
Using the second inequality in Lemma~\ref{nowy2}, we obtain that
\[
\biggl(\bigl|X\bigl(t^{l+1}_i\bigr)-X\bigl(t^l_i
\bigr)\bigr|-2^{-l+\bar{m}}\frac{c}{3} \biggr)_{+}\leq
b_{l,\bar{m}} \bigl[\varphi^{-1} \bigl(a_{l,\bar
{m}}^{-1}
\Delta\bigl(t^{l+1}_i,t^l_i\bigr)
\bigr) \bigr]^{\sfrac{1}{q}}.
\]
Now observe that $|t_i-t_{i-1}|\geq r^{-\vfrac{\bar{m}+1}{q}}S$ for
$i\in J_m$, $m\gs0$.
Therefore,
%
\begin{eqnarray}\label{weights1}
\sum^N_{l=\bar{m}+1} b_{l,\bar{m}}& =& \sum
^N_{l=\bar{m}+1} \biggl(\frac{c}{6}2^{l-\bar{m}}
\biggr)^{\vfrac{q-1}{q}} \bigl(C Lr^{-l}S^q
\bigr)^{\sfrac{1}{q}}
\nonumber
\\
&\leq& \bigl(6^{1-q}CL \bigr)^{\sfrac{1}{q}}\sum
^{\infty}_{l=\bar{m}+1} \bigl(2^{\vfrac{q-1}{q}}r^{-\sfrac{1}{q}}
\bigr)^{l-\bar{m}}c^{\vfrac
{q-1}{q}}r^{-\sfrac{\bar{m}}{q}}S \\
&\ls& M_2c^{\vfrac{q-1}{q}}|t_i-t_{i-1}|,
\nonumber
\end{eqnarray}
where $M_2 = M_2(\varphi,r,q)$ is defined by
\begin{eqnarray*}
 \bigl(6^{1-q}CL \bigr)^{\sfrac{1}{q}}r^{\sfrac{1}{q}}\sum
^{\infty
}_{l=\bar{m}+1} \bigl(2^{\vfrac{q-1}{q}}r^{-\sfrac{1}{q}}
\bigr)^{l-\bar{m}}&\leq& \bigl((12)^{1-q}CL\bigr)^{\sfrac{1}{q}}\sum
^{\infty}_{l'=0}\bigl(2^{\vfrac
{1-q}{q}}r^{-\sfrac{1}{q}}
\bigr)^{l'}
\\
&=&\bigl((12)^{1-q}CL\bigr)^{\sfrac{1}{q}}\bigl(1-2^{\vfrac{1-q}{q}}r^{-\sfrac
{1}{q}}
\bigr)^{-1}=:M_2.
\end{eqnarray*}
Consequently
\[
\sum^{\infty}_{m=0}\sum
_{i\in J_m}\sum^N_{l=\bar{m}+1}\sum
_{s\in
\{i-1,i\}}b_{l,\bar{m}} \ls2M_2c^{\vfrac{q-1}{q}}
\sum^n_{i=1}|t_i-t_{i-1}|=2M_2c^{\vfrac{q-1}{q}}S.
\]
Thus we can apply Lemma~\ref{convexf} for $\varphi(x^q)$ which is
convex above $C_{\varphi,q}$ and get
%
\begin{equation}
\label{duma1} V_{2}\leq2 M_2 c^{\vfrac{q-1}{q}}S \bigl[
\varphi^{-1}(\bar {V}_{2}+D_{\varphi,q})
\bigr]^{-1},
\end{equation}
where
\[
\bar{V}_{2}:= (2M_2)^{-1}\sum
^{\infty}_{m=0}\sum_{i\in J_m}
\sum^N_{l=\bar{m}+1}\sum
_{s\in\{i-1,i\}} \frac{\bar{b}_{l,\bar{m}}}{a_{l,\bar{m}}}\Delta\bigl(t^l_i,t^{l+1}_i
\bigr)
\]
and
\[
\bar{b}_{l,\bar{m}}=\bigl(c^{\vfrac{q-1}{q}}S\bigr)^{-1}b_{l,\bar{m}}=
\bigl(6^{-1}2^{-l+\bar{m}} \bigr)^{\vfrac{q-1}{q}}
\bigl(CLr^{-l} \bigr)^{\sfrac{1}{q}}.
\]
Observe that
\[
\bar{b}_{l,\bar{m}}\ls\bar{b}_{l,m_0}, a_{l,\bar{m}}\gs
a_{l,m_0},
\]
which implies $\bar{b}_{l,\bar{m}}/a_{l,\bar{m}}\ls\bar
{b}_{l,m_0}/a_{l,m_0}$.
Hence by Lemma~\ref{nowy1}, we have
%
\begin{equation}
\label{duma11} \bar{V}_{2}\leq\mathcal{V}_{2}:=M_2^{-1}
\sum^{\infty
}_{l=m_0+1}\frac{\bar{b}_{l,m_0}}{a_{l,m_0}} \sum
_{u\in T_{l+1}}\Delta\bigl(\pi_l(u),u\bigr).
\end{equation}
By the construction, $a_{l,m_0}\gs1$, furthermore by (\ref{ham0})
\begin{eqnarray*}
 M_2^{-1}\sum^{\infty}_{l=m_0+1}
\frac{\bar
{b}_{l,m_0}}{a_{l,m_0}}|T_{l+1}|&\ls& M_2^{-1}\sum
^{\infty}_{l=m_0+1} \frac{ (6^{-1}2^{-l+m_0} )^{\vfrac{q-1}{q}}
(CLr^{-l} )^{\sfrac{1}{q}}}{\varphi (\vafrac
{2^{-l+m_0}r^lc}{6CLS^q} )}
\bigl(r^{\vfrac{l+1}{q}}+1\bigr)
\\
&\ls&2M_2^{-1} \bigl((12)^{1-q}CL
\bigr)^{\sfrac{1}{q}} \sum^{\infty}_{l'=0}
\frac{2^{l'\vfrac{1-q}{q}}}{\varphi(2^{-l'}r^{l'})}=:M_3,
\end{eqnarray*}
where we have used the fact that $(r^{\vfrac{l+1}{q}}+1)\ls2r^{\vfrac
{l+1}{q}}$ and the definition
of $m_0$, that is, $r^{-m_0-1}S^q<c/M_0$, $M_0=12CL$, together with the
monotonicity of $\varphi$.
Note that by Remark~\ref{poziomka}, for convex $\varphi$, that is,
$C_{\varphi}=0$,
%
\begin{equation}
\label{pazur} \sum^{\infty}_{l'=0}
\frac{2^{l'\vfrac{1-q}{q}}}{\varphi
(2^{-l'}r^{l'})}\ls \sum^{\infty}_{l'=0}
4^{l'}r^{-l'}=\bigl(1-4r^{-1}\bigr)^{-1}.
\end{equation}
For $\varphi$ which is convex for $x\gs C_{\varphi}$ basically the
same argument works
but $l'$ must be large enough to apply the convexity. Indeed, using
that $\psi(x)=\varphi(x+C_{\varphi})-\varphi(C_\varphi)$ is convex
and $\psi(0)=0$ we deduce
$\psi(2^{-l'}r^{l'}x)\gs2^{-l'}r^{l'}\psi(x)$ for $x\gs0$ and thus
for all $x\gs0$,
%
\begin{equation}
\label{pantera} \varphi\bigl(2^{-l'}r^{l'}x+C_{\varphi}
\bigr)\gs2^{-l'}r^{l'}\bigl(\varphi (x+C_{\varphi})-
\varphi(C_{\varphi})\bigr)+\varphi(C_{\varphi}).
\end{equation}
Now choosing a suitable $x$ one can get a bound similar to (\ref{pazur})
yet for general $\varphi$. Note that in this case the bounding
constant may depend on $\varphi$. It proves that $M_3<\infty$.
Finally, by (\ref{duma1}), (\ref{duma11}) and Lemma~\ref{convexf1}
we get
%
\begin{equation}
\label{duma111} V_{2}\leq2 M_2 \max\{M_3,1\}
c^{\vfrac{q-1}{q}}S \bigl[\varphi ^{-1}({\mathcal{V}_2}/M_3+D_{\varphi,q})
\bigr]^{\sfrac{1}{q}}.
\end{equation}
Clearly, by (\ref{cond1}) and the definition of $M_3$ we have $\mathbf
{E}{\mathcal{V}_2/M_3} \leq1$.

A similar argument can be used to bound increments in $W_2$. Namely
using the forth inequality in Lemma~\ref{nowy2} we get that for $m>m_0$
and $i\in J_m$
\[
\biggl(\bigl|X\bigl(t^{m+1}_{i}\bigr)-X\bigl(t^{m+1}_{i-1}
\bigr)\bigr|-\frac{c}{3} \biggr)_{+}\ls b_m \bigl[
\varphi^{-1} \bigl(a_m^{-1}\Delta
\bigl(t^{m+1}_i,t^{m+1}_{i-1}\bigr)
\bigr) \bigr]^{\sfrac{1}{q}}.
\]
Using that $r^{-\vfrac{m+1}{q}}S\leq|t_i-t_{i-1}|\leq r^{-\sfrac
{m}{q}}S$ we get
\[
b_m= \biggl(\frac{c}{6} \biggr)^{\vfrac{q-1}{q}}
\bigl(2CLr^{-m}S^q \bigr)^{\sfrac{1}{q}}\ls
M_4c^{\vfrac{q-1}{q}}|t_i-t_{i-1}|,
\]
where $M_4= (2\cdot6^{1-q}CLr^{-1} )^{\sfrac{1}{q}}$.
Therefore,
\[
\sum^{\infty}_{m=m_0+1}\sum
_{i\in J_m}b_m\ls M_4 c^{\vfrac
{q-1}{q}}
\sum^n_{i=1}|t_i-t_{i-1}|=
M_4 c^{\vfrac{q-1}{q}}S,
\]
and thus using Lemma~\ref{convexf} for $\varphi(x^q)$ we get
%
\begin{equation}
\label{duma2} W_2\leq M_4 c^{\vfrac{q-1}{q}}S \bigl[
\varphi^{-1}(\bar {W}_2+D_{\varphi,q})
\bigr]^{\sfrac{1}{q}},
\end{equation}
where
\[
\bar{W}_2:=M_4^{-1}\sum
^{\infty}_{m=m_0+1}\sum_{i\in J_m}
\frac
{\bar{b}_m}{a_m}\Delta\bigl(t^{m+1}_i,t^{m+1}_{i-1}
\bigr)
\]
and $\bar{b}_m=(c^{\vfrac{q-1}{q}}S)^{-1}b_m=(2\cdot
6^{1-q}CLr^{-m})^{\sfrac{1}{q}}$. By (\ref{bakcyl3})
we have $d(t^{m+1}_i,t^{m+1}_{i-1})\ls2r{-m}S^q$ and thus using the
definition of the set $I_{m+1}(u)$ for each $m>m_0$ and $u\in T_{m+1}$
\[
\bar{W}_2\leq\mathcal{W}_2 = M_4^{-1}
\sum^{\infty}_{m=m_0+1} \frac
{\bar{b}_m}{a_m} \sum
_{u\in T_{m+1}}\sum_{v\in I_{m+1}(u)}
\Delta(u,v).
\]
Note that by (\ref{ham0}), (\ref{ham2})
\[
\sum_{u\in T_{m+1}}\bigl|I_{m+1}(u)\bigr|
\ls2^{-1}B(r,q) \bigl(r^{\vfrac
{m+1}{q}}+1\bigr)\ls B(r,q)r^{\vfrac{m+1}{q}}.
\]
Hence
\begin{eqnarray*}
&& M_4^{-1}B(r,q)\sum^{\infty}_{m=m_0+1}
\frac{\bar
{b}_m}{a_m}r^{\vfrac{m+1}{q}} M_4^{-1}B(r,q)\sum
^{\infty
}_{m=m_0+1}\frac{ (2\cdot6^{1-q}CLr  )^{\sfrac
{1}{q}}}{\varphi (\vafrac{r^mc}{12CLS^q} )}
\\
&&\quad \ls M_4^{-1}B(r,q) \bigl(2\cdot6^{1-q}CLr
\bigr)^{\sfrac{1}{q}}\sum^{\infty}_{m'=0}\bigl(
\varphi\bigl(r^{m'}\bigr)\bigr)^{-1}=:M_5,
\end{eqnarray*}
where in the last line we used that $r^{-m_0-1}S^q<c/M_0$, $M_0=12CL$.
The same argument as for $M_3$
proves that $M_5<\infty$. Note that in the case of convex $\varphi$
we can easily bound
$\sum^{\infty}_{m'=0}(\varphi(r^{m'}))^{-1}$ by $(1-r^{-1})^{-1}$.
By (\ref{duma2}) and Lemma~\ref{convexf1}, we get
%
\begin{equation}
\label{duma21} W_2\leq M_4 \max\{M_5,1\}
c^{\vfrac{q-1}{q}}S \bigl[\varphi ^{-1}(\mathcal{W}_2/M_5+D_{\varphi,q})
\bigr]^{\sfrac{1}{q}}.
\end{equation}
Obviously $\E\mathcal{W}_2/ M_5 \leq1$, consequently by
(\ref{duma111}), (\ref{duma21}) and Jensen's inequality we obtain the
desired result.
\end{pf}
Now we are ready to finish the proof of Theorem~\ref{Fund}.
\begin{pf*}{Proof of Theorem~\ref{Fund}}
Note that for fixed $q$ and $\varphi$ we may minimize constants
$K_1(\varphi,r,q)$ and $K_2(\varphi,r,q)$ appearing Lemmas \ref
{drugi}, \ref{trzeci} with respect to $r\geq4$.
It is clear from our discussion about the finitness of $M_3,M_5$ that
one can set $r=4$ in the case of convex $\varphi$. If $C_{\varphi}>0$
the choice of $r\gs4$ may be of
meaning as we have explained in (\ref{pantera}). Such minimal
constants depend only on $\varphi$ and $q$, and we will denote them by
$K_1(\varphi,q)$ and $K_2(\varphi,q)$ respectively.
Now it suffices to use Lemma~\ref{pierwszy}, then universal bounds
given in Lemmas \ref{drugi}, \ref{trzeci}
and finally let $N\rightarrow\infty$. Recall that by the construction
variables $Z_1$ and $Z_2$ of Lemmas \ref{drugi}, \ref{trzeci} do not
depend on $N$ and $\lim_{N\rightarrow\infty}d(t, \pi_{N+1}(t))=0$
for any $t\in T$. From condition (\ref{cond1}), for a given partition
$\Pi_n=\{t_0,t_1,\dots,t_n\}$ we get $R_N:=\sum^n_{i=1}|X(t_i)-X(t^{N+1}_i)| \rightarrow0$ in probability as $N
\uparrow+\infty$. Taking subsequence $N_k$ such that $R_{N_k}
\rightarrow0$ almost surely, we get the universal bound for the sum
$\sum^n_{i=1} (|X(t_i)-X(t_{i-1})|-c )_{+}$. Since $\Pi_n$
was arbitrary we get the result for $\TV^c (X,S )$.
\end{pf*}

\subsection{Application to the fractional Brownian motion}\label{sect3}

Let $W_H(t)$, $t\geq0$, be a fractional Brownian motion of the Hurst
parameter $H\in(0,1)$, that is, a centered Gaussian process which has the
following covariance function
%
\begin{equation}
\label{fBm_def} \mathbf{E} \bigl(W_H(s)W_H(t) \bigr) =
\tfrac{1}{2} \bigl(s^{2H} + t^{2H} - |s-t|^{2H}
\bigr).
\end{equation}
Let us consider $T=[0,S]$ with distance $d(s,t)=|t-s|^{H}$.
From (\ref{fBm_def}), it follows that $W_H(t) - W_H(s) \sim{\mathcal
N} (0, |t-s|^{2H}  )$ and thus, for some constant $C(H)$,
\[
\mathbf{E}\varphi_2 \biggl(\frac{|W_H(t)-W_H(s)|}{C(H)|t-s|^H} \biggr)\leq1,\qquad
\mbox{for } s,t\in T, s\neq t.
\]
Consequently, all assumptions of Corollary~\ref{cord2} are satisfied
with $p=2$, $q=H$ and we get
the following corollary.

\begin{coro}\label{las}
For any fractional Brownian motion $W_H(t)$, $t\in T$, the following
inequality holds
\[
\mathbf{P} \bigl(\TV^c(W_H,S)\geq c^{\vfrac{H-1}{H}}S(A_H+B_Hu)
\bigr)\leq C_H\exp \bigl(-u^{2H} \bigr),\qquad  \mbox{for } u>0,
\]
where $A_H,B_H,C_H$ are universal constants and $C_H=1$ for $H\gs1/2$.
\end{coro}

Note that Corollary~\ref{las} implies that $\mathbf{E}
\TV^c(W_H,S)\leq K_H c^{\vfrac{H-1}{H}}S$, where $K_H<\infty$.
On the other hand $c^{\vfrac{H-1}{H}}S$ is also the proper lower bound
for $\mathbf{E} \TV^c(W_H,S)$ when $S c^{-1/H} $ is not too small.
Indeed, let us consider the partition
$\Pi=\{0\leq t_0< t_1 <\cdots< t_N \leq S\}$ given by $t_i = i c^{1/H}$,
$i=0,1,2,\ldots, N=\lfloor S c^{-1/H} \rfloor$.
We have
\[
\TV^c(W_H,S)\geq\sum^N_{i=1}
\bigl(\bigl|W_H(t_i)-W_H(t_{i-1})\bigr|-c
\bigr)_{+}.
\]
Clearly, for $S c^{-1/H} \geq2$, $N > S c^{-1/H} -1 \geq S c^{-1/H}/2$
and $\mathbf{E}(|W_H(t_i)-W_H(t_{i-1})|-c)_{+} \geq k_H c$ for some
positive constant $k_H$.
It proves that when $S c^{-1/H} \geq2$, $c^{\vfrac{H-1}{H}}S$ is
comparable with $\mathbf{E} \TV^c(W_H,S)$ up to a constant depending
only on $H$.
Therefore, we have another formulation of Corollary~\ref{las}.

\begin{coro}
Assume that $S c^{-1/H} \geq2$. For any fractional Brownian motion
$W_H(t)$, $t\in T$, the following inequality holds
\[
\mathbf{P} \bigl(\TV^c(W_H,S)\geq\mathbf{E}
\TV^c(W_H,S) (\bar {A}_H+
\bar{B}_Hu) \bigr)\leq\bar{C}_H\exp\bigl(-u^{2H}
\bigr), \qquad \mbox{for } u>0,
\]
where $\bar{A}_H,\bar{B}_H,\bar{C}_H <\infty$ are universal
constants. Moreover $\bar{C}_H=1$ for $H\geq1/2$.
\end{coro}

\section{Application to the standard Brownian motion and
diffusions}\label{sect4}

For a standard Brownian motion $W=W_{1/2}$, which is the only fractional
Brownian motion with independent increments one may, using this property,
strengthen the results obtained for general fBm and obtain Gaussian
concentration of $\TV^{c} (W,S )$. The generalization of this
result for diffusions with moderate growth, driven by $W$, is also possible.

Let us assume that $X_{t}$, $t\geq0$, is a one-dimensional diffusion satisfying
%
\begin{equation}\label{eq:sde}
X (t )=x_{0}+\int_{0}^{t}\mu
\bigl(s,X (s ) \bigr)\,\mathrm{d}s+\int_{0}^{t}\sigma
\bigl(s,X (s ) \bigr)\,\mathrm{d}W (s ).
\end{equation}
We assume that $\sigma\dvtx  [0;+\infty )\times\mathbb
{R}\rightarrow [-R;R ]$
is measurable and bounded (i.e., $0<R<+\infty$) and $\mu\dvtx
[0;+\infty )\times\mathbb{R}\rightarrow\mathbb{R}$
is measurable and satisfying the following linear growth condition: there
exists $C,D \geq0$ such that for all $t\geq0$
%
\begin{equation}\label
{eq:lipschitz}
\bigl\llvert \mu (t,x )\bigr\rrvert \leq C+D\llvert x\rrvert .
\end{equation}
We will also need the natural assumption that $X$ is a Markov process.
With this assumption, we have the following theorem.

\begin{theorem}\label{diffusions}
For $X$ being a Markov process satisfying (\ref{eq:sde}) with
$\mu$ and $\sigma$ as above and $\lambda\geq0$ one has
\begin{eqnarray*}
\mathbf{E}\exp \bigl(\lambda \TV^{c} (X,S ) \bigr) & \leq & 2\exp
\bigl(\lambda^{2}S\alpha_{R}+\lambda Sc^{-1}\beta
_{R}+\lambda\gamma_{x_{0},C,D,S} \bigr)
\\
& &{} \times \bigl(1+8\lambda\eta_{D,R,S}\exp \bigl(\lambda^{2}
\eta _{D,R,S}^{2} \bigr) \bigr),
\end{eqnarray*}
where $\gamma_{x_{0},C,D,S}= (C+D\llvert x_{0}\rrvert  )S\mathrm{e}^{DS}$,
$\delta_{D,S}=DS\mathrm{e}^{DS}$ and $\eta_{D,R,S}=\delta_{D,S}R\sqrt{S/2}$.
In particular, when $D=0$ we get
\[
\mathbf{E}\exp \bigl(\lambda \TV^{c} (X,S ) \bigr)\leq 2\exp \bigl(
\lambda^{2}S\alpha_{R}+\lambda S \bigl(c^{-1}
\beta _{R}+C \bigr) \bigr)
\]
and for the standard Brownian motion $X=W$ we get
%
\begin{equation}
\label{wiener} \mathbf{E}\exp \bigl(\lambda \TV^{c} (W,S ) \bigr)\leq 2
\exp \bigl(\lambda^{2}S\alpha+\lambda Sc^{-1}\beta \bigr),
\end{equation}
where $\alpha,\beta$ are universal constants.
\end{theorem}

\begin{pf}
Let us define
\[
M (t ):=\int_{0}^{t}\mu \bigl(s,X (s )
\bigr)\,\mathrm{d}s,\qquad Y (t ):=\int_{0}^{t}\sigma \bigl(s,X (s )
\bigr)\,\mathrm{d}W (s )
\]
and $Y^{*}=\sup_{0\leq s\leq S}\llvert Y (s )\rrvert $. We
have $X (t )=x_{0}+M (t )+Y (t )$, and
due to (\ref{eq:lipschitz}) we estimate
%
\begin{eqnarray}\label{eq:estim}
\bigl\llvert M (t )\bigr\rrvert & \leq& \int_{0}^{t}
\bigl\llvert \mu \bigl(s,X (s ) \bigr)\bigr\rrvert \,\mathrm{d}s\leq\int_{0}^{t}C+D
\bigl\llvert X (s )\bigr\rrvert\, \mathrm{d}s
\nonumber
\\
& \leq& \int_{0}^{t}C+D\llvert x_{0}
\rrvert +D\bigl\llvert M (s )\bigr\rrvert +DY^{*}\,\mathrm{d}s
\\
& \leq& \bigl(C+D\llvert x_{0}\rrvert +DY^{*} \bigr)S+D
\int_{0}^{t}\bigl\llvert M (s )\bigr\rrvert \,\mathrm{d}s.\nonumber
\end{eqnarray}
Hence, from Gronwall's lemma (cf. Revuz and Yor \cite{RY}, Appendix
\S1), we get
%
\begin{equation}\label{eq:gronwall}
\bigl\llvert M (t )\bigr\rrvert \leq \bigl(C+D\llvert x_{0}\rrvert
+DY^{*} \bigr)S\mathrm{e}^{Dt}.
\end{equation}
Notice that due to (\ref{eq:gronwall}) $M$ is adapted, absolute
continuous process with locally bounded total variation. Indeed, repeating
estimates (\ref{eq:estim}) and using (\ref{eq:gronwall}) we get
%
\begin{eqnarray}\label{eq:tms}
\TV (M,S ) & \leq& \int_{0}^{S}\bigl\llvert \mu
\bigl(s,X (s ) \bigr)\bigr\rrvert \,\mathrm{d}s
\nonumber
\\
& \leq& \bigl(C+D\llvert x_{0}\rrvert +DY^{*} \bigr)S+D
\int_{0}^{S}\bigl\llvert M (t )\bigr\rrvert \,\mathrm{d}s
\nonumber
\\[-8pt]\\[-8pt]
& \leq& \bigl(C+D\llvert x_{0}\rrvert +DY^{*} \bigr)S+D
\bigl(C+D\llvert x_{0}\rrvert +DY^{*} \bigr)S\int
_{0}^{S}\mathrm{e}^{Dt}\,\mathrm{d}s
\nonumber
\\
& = & \bigl(C+D\llvert x_{0}\rrvert \bigr)S\mathrm{e}^{DS}+DS\mathrm{e}^{DS}Y^{*}.\nonumber
\end{eqnarray}
($\TV=\TV^{0}$ denotes here the total variation.)

By \L ochowski and Mi\l o\'{s} \cite{LM}, Fact~17, we
have
%
\begin{equation}
\TV^{c} (X,S )\leq \TV (M,S )+\TV^{c} (Y,S ).\label{eq:tvcdecomp}
\end{equation}
Now we will investigate $\TV^{c} (Y,S )$.

First, let us prove
that $Y$ satisfies condition (\ref{lpchwyt}) with $\varphi= \varphi_2$ and
$d (s,t )=\llvert s-t\rrvert ^{1/2}$.
Indeed, let us fix $0\leq s<t\leq S$ and consider the following martingale
$Z (u ):=Y (s+u )-Y (s ),u\in
[0;t-s ]$.
We have
\[
Z (u )=\int_{s}^{s+u}\sigma \bigl(\tau,X (\tau )
\bigr)\,\mathrm{d}W (\tau )
\]
and
\[
\langle Z \rangle (u )=\int_{s}^{s+u}\sigma
\bigl(\tau,X (\tau ) \bigr)^{2}\,\mathrm{d}\tau\leq R^{2} (t-s ).
\]
Hence, by Bernstein's inequality (cf. Revuz and Yor \cite{RY},
Chapter IV, Exercise 3.16),
we have
%
\begin{eqnarray}\label{eq:bernstein}
 \mathbf{P} \bigl(\bigl\llvert Y (t )-Y (s )\bigr\rrvert \geq x \bigr)&\leq&2
\mathbf{P} \Bigl(\sup_{u\in [0;t-s
]}Z (u )\geq x \Bigr)
\nonumber
\\
& =& 2\mathbf{P} \Bigl(\sup_{u\in [0;t-s ]}Z (u )\geq x, \langle Z
\rangle (t-s )\leq R^{2} (t-s ) \Bigr)
\\
& \leq&2\exp \bigl(-x^{2}/ \bigl(2R^{2} (t-s ) \bigr)
\bigr).\nonumber
\end{eqnarray}
From (\ref{eq:bernstein}), we immediately get that $Y$ satisfies
condition (\ref{lpchwyt}) for $\varphi= \varphi_2$ and $d (s,t
)=\llvert s-t\rrvert ^{1/2}$.
Hence, from Corollary~\ref{cord2} we obtain the following bound on the tails
of $\TV^{c} (Y,S )$:
%
\begin{equation}\label{eq:corrollary1}
\mathbf{P} \bigl(\TV^{c} (Y,S )\geq c^{-1}S (A+Bu ) \bigr)
\leq \mathrm{e}^{-u},
\end{equation}
where $A=A (R )$ and $B=B (R )$ depend on $R$
only. Notice that for $\delta>0$ applying Bernstein's inequality to
$Y^{*}$ we get
$\mathbf{P} (Y^{*}\geq x )\leq2\exp
(-x^2/(2R^2S) )$
and using integration by parts we have
%
\begin{equation}\label{eq:bernstein1}
\mathbf{E}\exp \bigl(\delta Y^{*} \bigr)\leq1+2\delta\int
_{0}^{\infty}\mathrm{e}^{\delta y}\mathrm{e}^{-y^{2}/(2R^{2}S)}\,\mathrm{d}y\leq1+8
\delta R\sqrt {S/2}\mathrm{e}^{\delta^{2}R^{2}S/2}.
\end{equation}

Now, we will strengthen estimate (\ref{eq:corrollary1}) using the
Markov property of $X$. First, using (\ref{eq:corrollary1}) and
integration by parts we have
%
\begin{equation}\label{eq:expmoment}
\mathbf{E}\exp \bigl(\lambda \bigl[\TV^{c} (Y,S )-c^{-1}SA
\bigr] \bigr)\leq\frac{1}{1-\lambda SB/c}
\end{equation}
for $\lambda<c (SB )^{-1}$. Let now $S=S_{1}+S_{2}$, where
$S_{1},S_{2}>0$. Using the inequality $\TV^{c} (Y,S )\leq
\TV^{c} (Y,S_{1} )+c+\TV^{c} (Y, [S_{1},S
] )$,
which follows easily from the estimate:
\[
\bigl(\bigl\llvert Y (t )-Y (u )\bigr\rrvert -c \bigr)_{+}\leq
\bigl(\bigl\llvert Y (t )-Y (S_{1} )\bigr\rrvert -c
\bigr)_{+}+ \bigl(\bigl\llvert Y (S_{1} )-Y (u )\bigr
\rrvert -c \bigr)_{+}+c
\]
for $0\leq t<S_{1}<u\leq S$, and then the Markov property of $X$ we
get
%
\begin{eqnarray}\label{eq:expmoment2}
&& \mathbf{E} \exp \bigl(\lambda \bigl[\TV^{c} (Y,S
)-c^{-1}SA \bigr] \bigr)
\nonumber
\\
&&\quad  \leq\mathbf{E} \exp \bigl(\lambda \TV^{c} (Y,S_{1} )+
\lambda c+\lambda \TV^{c} \bigl(Y, [S_{1},S ] \bigr)-\lambda
c^{-1}SA \bigr)
\nonumber
\\[-8pt]\\[-8pt]
&&\quad  = \mathrm{e}^{\lambda c}\mathbf{E} \bigl(\mathrm{e}^{ \lambda \TV^{c}
(Y,S_{1} )-\lambda c^{-1}S_{1}A} \E \bigl[\mathrm{e}^{ \lambda
\TV^{c} (Y, [S_{1},S ] )-\lambda c^{-1}S_{2}A}
|X (S_{1} ) \bigr] \bigr)
\nonumber
\\
&&\quad  \leq \mathrm{e}^{\lambda c}\frac{1}{1-\lambda S_{1}B/c}\frac{1}{1-\lambda
S_{2}B/c}.\nonumber
\end{eqnarray}
The last inequality follows by (\ref{eq:expmoment}), since the right-hand
side of (\ref{eq:expmoment}) does not depend on $x_{0}$, and
using the Markov property in similar way we have the universal estimate
for the conditional expectation
\[
\E \bigl(\exp \bigl\{ \lambda \TV^{c} \bigl(Y, [S_{1},S ]
\bigr)-\lambda c^{-1}S_{2}A \bigr\} |X (S_{1}
)=x_{1} \bigr)\leq\frac{1}{1-\lambda S_{2}B/c}
\]
(note that the length of interval $ [S_{1},S ]$ is $S_{2}$).
Notice now that from (\ref{eq:expmoment2}) it follows that $\mathbf
{E} \exp (\lambda [\TV^{c} (Y,S )-c^{-1}S\bar
{A} ] )<+\infty$
for $\lambda<\min \{ c (S_{1}B )^{-1},c
(S_{2}B )^{-1} \} $.
Let us fix integer $n\geq1$. Iterating (\ref{eq:expmoment2}) we
obtain
%
\begin{equation}\label{eq:estimn}
\mathbf{E} \exp \bigl(\lambda \bigl[\TV^{c} (Y,S )-c^{-1}SA
\bigr] \bigr)\leq \mathrm{e}^{\lambda c (n-1 )} \biggl(\frac{1}{1-\lambda SB (cn )^{-1}}
\biggr)^{n}
\end{equation}
for $\lambda<cn (SB )^{-1}$, which gives that $\mathbf{E}
\exp (\lambda [\TV^{c} (Y,S )-c^{-1}SA
] )<+\infty$
for any $\lambda\in\mathbb{R}$. Now, let us fix $\lambda>0 $
and set $n= \lceil2\lambda SBc^{-1} \rceil$. Using (\ref
{eq:estimn}),
we get
\begin{eqnarray*}
\mathbf{E} \exp \bigl(\lambda \bigl[\TV^{c} (Y,S )-c^{-1}SA
\bigr] \bigr) & \leq& \mathrm{e}^{\lambda c (n-1
)}2^{n}
\\
& \leq& 2\exp \bigl(2\lambda^{2}SB+2 (\ln2 )\lambda
SBc^{-1} \bigr)
\end{eqnarray*}
and thus
%
\begin{eqnarray}\label{eq:estimtvy}
\mathbf{E} \exp \bigl(\lambda \TV^{c} (Y,S ) \bigr) & \leq & 2\exp
\bigl(2\lambda^{2}SB+\lambda Sc^{-1} \bigl(A+2 (\ln 2 )B
\bigr) \bigr)
\nonumber
\\[-8pt]\\[-8pt]
& = & 2\exp \bigl(\lambda^{2}S\alpha_{R}+\lambda
Sc^{-1}\beta _{R} \bigr),\nonumber
\end{eqnarray}
where $\alpha_{R}=2B=2B (R )$ and $\beta_{R}=A+2 (\ln
2 )B=A (R )+2 (\ln2 )B (R )$.
Now, from (\ref{eq:tvcdecomp}), (\ref{eq:tms}) and (\ref{eq:estimtvy})
we get
\begin{eqnarray*}
\mathbf{E} \exp \bigl(\lambda \TV^{c} (X,S ) \bigr) & \leq & E\exp
\bigl(\lambda \TV (M,S )+\lambda \TV^{c} (Y,S ) \bigr)
\\
& \leq& 2\exp \bigl(\lambda^{2}S\alpha_{R}+\lambda
Sc^{-1}\beta _{R}+\lambda\gamma_{x_{0},C,D,S} \bigr)\E
\exp \bigl(\lambda\delta _{D,S}Y^{*} \bigr),
\end{eqnarray*}
where $\gamma_{x_{0},C,D,S}= (C+D\llvert x_{0}\rrvert  )S\mathrm{e}^{DS}$,
$\delta_{D,S}=DS\mathrm{e}^{DS}$.
Finally, using (\ref{eq:bernstein1}) with $\delta= \lambda\delta
_{D,S}$ we get
\begin{eqnarray*}
\mathbf{E} \exp \bigl(\lambda \TV^{c} (X,S ) \bigr) & \leq & 2\exp
\bigl(\lambda^{2}S\alpha_{R}+\lambda Sc^{-1}\beta
_{R}+\lambda\gamma_{x_{0},C,D,S} \bigr)
\\
& &{}\times \bigl(1+8\lambda\eta_{D,R,S}\exp \bigl(\lambda^{2}
\eta _{D,R,S}^{2} \bigr) \bigr),
\end{eqnarray*}
where $\eta_{D,R,S}=\delta_{D,S}R\sqrt{S/2}$.
\end{pf}

\begin{rema}
Let us notice that the condition that $\sigma$ is bounded
is essential for obtaining the Gaussian concentration of $\TV^{c}
(X,S )$.
To see this it is enough to consider the equation $dX (t
)=2^{-1}X(t)\,\mathrm{d}t+X (t )\,\mathrm{d}W (t )$
with the starting condition $X (0 )=1$. Notice that
$\TV^{c} (X,S )\geq (X (S )-X(0)-c )_{+}$
and that $ (X (S )-X(0)-c )_{+}= (\exp
W (S )-1-c )_{+}$
does not reveal the Gaussian concentration.
\end{rema}

\begin{rema}
Notice that for the standard Brownian motion $X=W$ and $Sc^{-2} \geq
2$, $Sc^{-1}$ is comparable up to a universal constant with $\mathbf
{E}\TV^{c} (W,S )$. Hence, from (\ref{wiener}) we obtain
that for $c>0$ such that $Sc^{-2} \geq2$, there exist universal
constants $ \bar{A}, \bar{B} <+\infty$
such that the Gaussian concentration holds
\begin{eqnarray*}
\mathbf{P}\bigl(\TV^c(W,S)\geq\bar{A} \mathbf{E} \TV^c(W,S)+
\bar{B}\sqrt {S}u \bigr)\leq\exp\bigl(-u^2\bigr),\qquad  \mbox{for } u\geq0.
\end{eqnarray*}
\end{rema}

\section{Existence of moment-generating functions of the truncated
variation of L\'{e}vy processes}\label{sect5}

In this section, we will deal with the existence of finite exponential
moments of the truncated variation of a L\'{e}vy process $X$. We will
state the
necessary and sufficient condition for the finiteness of
$
\mathbf{E}\exp (\alpha \TV^{c} (X,S ) )
$ in terms of the generating triplet of the process $X$ (cf.
Sato \cite{Sato:1999}, Chapter~2, Section~11). The
methodology used here is very
similar to the methodology
used in \L ochowski \cite{Loch1} for a Wiener process $W$,
where the existence of
$\mathbf{E}\exp (\alpha \TV^{c} (W,S ) )$
for any complex $\alpha$ was proved.

We start with the following lemma.

\begin{lema} \label{lema8}
Let $X$ be a L\'{e}vy process. For any $c>0$ and $\alpha>0$ one has
$\mathbf{E}\exp (\alpha \TV^{c} (X,S ) ) <
+\infty$
if and only if
\[
\mathbf{E}\exp \Bigl(\alpha\sup_{0\leq s\leq S}\bigl\llvert X (s )\bigr
\rrvert \Bigr)<+\infty.
\]
\end{lema}

\begin{pf}
The `only if' part follows from the inequality
\begin{eqnarray*}
\TV^{c} (X,S ) & \geq&\sup_{0\leq s\leq S}\max \bigl\{ \bigl
\llvert X (s )- X (0 )\bigr\rrvert -c,0 \bigr\}
\\
& =&\max\Bigl\{ \sup_{0\leq s\leq S}\bigl|X(s)\bigr |-c,0\Bigr\} \geq\sup
_{0\leq s\leq S}\bigl|X(s)\bigr|-c.
\end{eqnarray*}
To prove the opposite implication let us define $T_{0}^{c}=0$
and for $i=1,2,\dots$
\[
T_{i}^{c}= \inf \bigl\{ t>T_{i-1}^{c}\dvtx
\bigl\llvert X (t )-X \bigl(T_{i-1}^{c} \bigr)\bigr\rrvert
>c/2 \bigr\} \wedge \bigl(S+T_{i-1}^{c} \bigr).
\]
Observe that $T_{1}^{c}=\inf \{ t>0\dvtx \llvert X (t
)\rrvert >c/2 \} \wedge S\leq S$
and that $ (X (t ) )_{t\geq0}\,\displaystyle \mathop{=}^{d}\, (X
(t )-X (T_{1}^{c} ) )_{t\geq T_{1}^{c}}$,
where ``$\displaystyle \mathop{=}^{d}$'' denotes the equality of distributions. Now let us define
\[
X_{t}^{c}=\sum_{i=0}^{\infty}X
\bigl(T_{i}^{c} \bigr)I_{
[T_{i}^{c},T_{i+1}^{c} )} (t ).
\]
Since $\llVert  X^{c}-X\rrVert _{\infty}\leq c/2$, we have
%
\begin{equation}\label{eq:tvcxleqtvxc}
\TV^{c} (X,S )\leq \TV \bigl(X^{c},S \bigr)
\end{equation}
and since $X^{c}$ is piecewise constant with the first jump at
$T_{1}^{c}\leq S$, denoting $\Delta X^c (T_{1}^{c} ) =
X^c (T_{1}^{c} ) - X^c (T_{1}^{c}- )$
we have
%
\begin{eqnarray}\label
{eq:tvxcleqsum}
\TV \bigl(X^{c},S \bigr) & =&\bigl\llvert \Delta X^{c}
\bigl(T_{1}^{c} \bigr)\bigr\rrvert +\TV
\bigl(X^{c}, \bigl[T_{1}^{c},S \bigr]
\bigr)\nonumber
\\[-8pt]\\[-8pt]
& \leq&\sup_{0\leq s\leq T_{1}^{c}}\bigl\llvert X (s )\bigr\rrvert +\TV
\bigl(X^{c}, \bigl[T_{1}^{c},S \bigr] \bigr).
\nonumber
\end{eqnarray}

Let now $\delta\in (0;S )$ be such a small number that
%
\begin{eqnarray}\label{eq:deltacond1}
&&\mathbf{E} \Bigl[\exp \Bigl(\alpha\sup_{0\leq s\leq S}\bigl\llvert X
(s )\bigr\rrvert \Bigr);T_{1}^{c}\leq\delta \Bigr]
\nonumber
\\[-8pt]\\[-8pt]
&&\quad :=\mathbf{E} \Bigl[\exp \Bigl(\alpha\sup_{0\leq s\leq S}\bigl\llvert X
(s )\bigr\rrvert \Bigr)I_{ \{ T_{1}^{c}\leq\delta
\} } \Bigr]<1. \nonumber
\end{eqnarray}
Note that such a number exists, since we assume that $\mathbf{E}\exp
(\alpha\sup_{0\leq s\leq S}\llvert X (s )\rrvert  )<+\infty$
and from the c\`{a}dl\`{a}g property and stochastic continuity of $X$ it
follows that $\mathbf{P} (T_{1}^{c}\leq\delta )=\mathbf
{P} (\sup_{0\leq s\leq\delta}\llvert X (s )\rrvert >c/2 )\downarrow0$
as $\delta\downarrow0$.

Let us fix $M>0$. Note that on the set $ \{T_{1}^{c}>\delta
\}$ we have $\TV (X^{c},\delta ) = 0$, hence
%
\begin{eqnarray}\label{phi}
\mathbf{E}\exp \bigl(\alpha \TV \bigl(X^{c},\delta \bigr)\wedge M
\bigr) & =&\mathbf{E} \bigl[\exp \bigl(\alpha \TV \bigl(X^{c},\delta
\bigr)\wedge M \bigr);T_{1}^{c}\leq\delta \bigr]
\nonumber
\\
&&{}+ \mathbf{E} \bigl[\exp (0 \wedge M );T_{1}^{c}>\delta
\bigr]
\nonumber
\\[-8pt]\\[-8pt]
& =&\mathbf{E} \bigl[\exp \bigl(\alpha \TV \bigl(X^{c},\delta \bigr)
\wedge M \bigr);T_{1}^{c}\leq\delta \bigr]
\nonumber
\\
&&{}+\mathbf{P} \bigl(T_{1}^{c}>\delta \bigr).\nonumber
\end{eqnarray}
Now, applying (\ref{eq:tvxcleqsum}), the independence of the process
$X (t )-X (T_{1}^{c} ), t\geq T_{1}^{c}$,
and the two-dimensional r.v. $ (\sup_{0\leq s\leq T_{1}^{c}}\llvert X (s )\rrvert ,T_{1}^{c} )$
(to see this notice that $T_{1}^{c}$ is a stopping time and use the
strong Markov property of L\'{e}vy processes) and the equality of distributions
of $\TV (X^{c},s  )$ and $\TV (X^{c},
[T_{1}^{c};T_{1}^{c}+s ] )$
for any $s\geq0$, we have
\begin{eqnarray*}
&&\mathbf{E} \bigl[\exp \bigl(\alpha \TV \bigl(X^{c},\delta \bigr)\wedge
M \bigr);T_{1}^{c}\leq\delta \bigr]
\\
&&\quad \leq\mathbf{E} \Bigl[\exp \Bigl(\alpha \Bigl(\sup_{0\leq s\leq
T_{1}^{c}}\bigl
\llvert X (s )\bigr\rrvert +\TV \bigl(X^{c}; \bigl[T_{1}^{c};
\delta \bigr] \bigr) \Bigr)\wedge M \Bigr);T_{1}^{c}\leq
\delta \Bigr]
\\
&&\quad \leq\mathbf{E} \Bigl[\exp \Bigl(\alpha\sup_{0\leq s\leq
T_{1}^{c}}\bigl\llvert
X (s )\bigr\rrvert +\alpha \TV \bigl(X^{c}; \bigl[T_{1}^{c};
\delta+T_{1}^{c} \bigr] \bigr)\wedge M
\Bigr);T_{1}^{c}\leq\delta \Bigr]
\\
&&\quad =\mathbf{E} \Bigl[\exp \Bigl(\alpha\sup_{0\leq s\leq
T_{1}^{c}}\bigl\llvert X
(s )\bigr\rrvert \Bigr);T_{1}^{c}\leq\delta \Bigr]
\mathbf{E}\exp \bigl(\alpha \TV \bigl(X^{c},\delta \bigr)\wedge M
\bigr)+\mathbf{P} \bigl(T_{1}^{c}>\delta \bigr)
\\
&&\quad \leq\mathbf{E} \Bigl[\exp \Bigl(\alpha\sup_{0\leq s\leq S}\bigl\llvert
X (s )\bigr\rrvert \Bigr);T_{1}^{c}\leq\delta \Bigr]
\mathbf {E}\exp \bigl(\alpha \TV \bigl(X^{c},\delta \bigr)\wedge M
\bigr).
\end{eqnarray*}
By this and by (\ref{phi}), (\ref{eq:deltacond1}) we have\vspace*{1pt}
%
\begin{equation}\label{eq:est1}
\mathbf{E}\exp \bigl(\alpha \TV \bigl(X^{c},\delta \bigr)\wedge M
\bigr)\leq\frac{\mathbf{P} (T_{1}^{c}>\delta
)}{1-\mathbf{E} [\exp (\alpha\sup_{0\leq s\leq S}\llvert X (s )\rrvert  );T_{1}^{c}\leq\delta
]}.
\end{equation}
Using similar arguments as before (i.e., (\ref{eq:tvxcleqsum}),
independence of $X (t )-X (T_{1}^{c} ),t\geq T_{1}^{c}$,
and $ (\sup_{0\leq s\leq T_{1}^{c}}\llvert X (s )\rrvert ,T_{1}^{c} )$
and the equality of distributions
of $\TV (X^{c},s )$ and $\TV (X^{c},
[T_{1}^{c};T_{1}^{c}+s ] )$
for $s\geq0$) we obtain\vspace*{1pt}
\begin{eqnarray*}
&&\mathbf{E}\exp \bigl(\alpha \TV \bigl(X^{c},S \bigr)\wedge M \bigr)
\\
&&\quad \leq\mathbf{E} \Bigl[\exp \Bigl(\alpha\sup_{0\leq s\leq
T_{1}^{c}}\bigl\llvert
X (s )\bigr\rrvert +\alpha \TV \bigl(X^{c}; \bigl[T_{1}^{c};S+T_{1}^{c}
\bigr] \bigr)\wedge M \Bigr);T_{1}^{c}\leq \delta \Bigr]
\\
&&\qquad {}+\mathbf{E} \Bigl[\exp \Bigl(\alpha\sup_{0\leq s\leq
T_{1}^{c}}\bigl\llvert X
(s )\bigr\rrvert +\alpha \TV \bigl(X^{c}; \bigl[T_{1}^{c};S+T_{1}^{c}-
\delta \bigr] \bigr)\wedge M \Bigr);T_{1}^{c}>\delta \Bigr]
\\
&&\quad =\mathbf{E} \Bigl[\exp \Bigl(\alpha\sup_{0\leq s\leq
T_{1}^{c}}\bigl\llvert X
(s )\bigr\rrvert \Bigr);T_{1}^{c}\leq\delta \Bigr]
\mathbf{E} \bigl[\exp \bigl(\alpha \TV \bigl(X^{c},S \bigr)\wedge M
\bigr) \bigr]
\\
&&\qquad {}+\mathbf{E} \Bigl[\exp \Bigl(\alpha\sup_{0\leq s\leq
T_{1}^{c}}\bigl\llvert X
(s )\bigr\rrvert \Bigr);T_{1}^{c}>\delta \Bigr]\mathbf{E}
\bigl[\exp \bigl(\alpha \TV \bigl(X^{c},S-\delta \bigr)\wedge M \bigr)
\bigr]
\\
&&\quad \leq\mathbf{E} \Bigl[\exp \Bigl(\alpha\sup_{0\leq s\leq S}\bigl\llvert
X (s )\bigr\rrvert \Bigr);T_{1}^{c}\leq\delta \Bigr]
\mathbf {E} \bigl[\exp \bigl(\alpha \TV \bigl(X^{c},S \bigr)\wedge M
\bigr) \bigr]
\\
&&\qquad {}+\mathbf{E}\exp \Bigl(\alpha\sup_{0\leq s\leq S}\bigl\llvert X (s )
\bigr\rrvert \Bigr)\mathbf{E} \bigl[\exp \bigl(\alpha \TV \bigl(X^{c},S-
\delta \bigr)\wedge M \bigr) \bigr].
\end{eqnarray*}
From this, we have\vspace*{1pt}
\begin{eqnarray*}
&&\mathbf{E}\exp \bigl(\alpha \TV \bigl(X^{c},S \bigr)\wedge M \bigr)
\\
&&\quad \leq\frac{\mathbf{E}\exp (\alpha\sup_{0\leq s\leq S}\llvert X (s )\rrvert  )}{1-\mathbf{E} [\exp
(\alpha\sup_{0\leq s\leq S}\llvert X (s )\rrvert
);T_{1}^{c}\leq\delta ]}\mathbf{E}\exp \bigl(\alpha \TV \bigl(X^{c},S-
\delta \bigr)\wedge M \bigr).
\end{eqnarray*}
Similarly, if $S-2\delta>0$\vspace*{1pt}
\begin{eqnarray*}
&&\mathbf{E}\exp \bigl(\alpha \TV \bigl(X^{c},S-\delta \bigr)\wedge M
\bigr)
\\
&&\quad \leq\frac{\mathbf{E}\exp (\alpha\sup_{0\leq s\leq S-\delta
}\llvert X (s )\rrvert  )}{1-\mathbf{E} [\exp
(\alpha\sup_{0\leq s\leq S-\delta}\llvert X (s
)\rrvert  );T_{1}^{c}\leq\delta ]}\mathbf{E}\exp \bigl(\alpha \TV \bigl(X^{c},S-2
\delta \bigr)\wedge M \bigr)
\\
&&\quad \leq\frac{\mathbf{E}\exp (\alpha\sup_{0\leq s\leq S}\llvert X (s )\rrvert  )}{1-\mathbf{E} [\exp
(\alpha\sup_{0\leq s\leq S}\llvert X (s )\rrvert
);T_{1}^{c}\leq\delta ]}\mathbf{E}\exp \bigl(\alpha \TV \bigl(X^{c},S-2
\delta \bigr)\wedge M \bigr).
\end{eqnarray*}
Iterating and putting together the above inequalities, we finally
obtain\vspace*{1pt}
%
\begin{eqnarray}\label{eq:est3}
\mathbf{E}\exp \bigl(\alpha \TV \bigl(X^{c},S \bigr)\wedge M \bigr) &
\leq &\biggl(\frac{\mathbf{E}\exp (\alpha\sup_{0\leq s\leq
S}\llvert X (s )\rrvert  )}{1-\mathbf{E} [\exp
(\alpha\sup_{0\leq s\leq S}\llvert X (s )\rrvert  );T_{1}^{c} \leq\delta ]} \biggr)^{ \lfloor
S/\delta \rfloor}
\nonumber
\\[-8pt]\\[-8pt]
&&{} \times\mathbf{E}\exp \bigl(\alpha \TV \bigl(X^{c},\delta \bigr)
\wedge M \bigr).\nonumber
\end{eqnarray}
By (\ref{eq:est1}) and (\ref{eq:est3}), and letting $M\rightarrow
\infty$
we get $\mathbf{E}\exp (\alpha \TV (X^{c},S )
)<+\infty$.
Finally, from (\ref{eq:tvcxleqtvxc}) we get
\[
\mathbf{E}\exp \bigl(\alpha \TV^{c} (X,S ) \bigr)<+\infty.
\]
\upqed
\end{pf}
Now let $ (A,\nu,\gamma )$ be the generating triplet of the
process $X$. By Sato \cite{Sato:1999}, Theorem~28.15, we have
\[
\mathbf{E}\exp \Bigl(\alpha\sup_{0\leq s\leq S}\bigl\llvert X (s )\bigr
\rrvert \Bigr)<+\infty
\]
if and only if
\begin{eqnarray*}
\mathbf{E}\exp \bigl(\alpha\bigl\llvert X (1 )\bigr\rrvert \bigr)<+\infty
\end{eqnarray*}
which, by Sato \cite{Sato:1999}, Corollary~25.8,
is equivalent with
%
\begin{equation}
\label{levy_cond} \int_{\llvert x\rrvert >1}\mathrm{e}^{\alpha\llvert x\rrvert }\nu (\mathrm{d}x )<+\infty.
\end{equation}
From equivalence of these conditions and Lemma~\ref{lema8} we obtain
the following theorem.

\begin{theorem}\label{Levy}
Let $ (A,\nu,\gamma )$ be the generating triplet of the L\'{e}vy
process $X$. For any $\alpha>0$ we have
\[
\mathbf{E}\exp \bigl(\alpha \TV^{c} (X,S ) \bigr) < +\infty
\]
if and only if
\[
\int_{\llvert x\rrvert >1}\mathrm{e}^{\alpha\llvert x\rrvert }\nu (\mathrm{d}x )<+\infty.
\]
\end{theorem}

Theorem~\ref{Levy} may be applied in situations, when the process $X$
satisfies condition (\ref{levy_cond}) with some $\alpha>0$ but it is
neither Brownian motion nor finite variation process. This holds, for example,
for tempered stable process, that is, processes with the L\'{e}vy
measure given by
\[
\nu(\mathrm{d}x) = \frac{c_p}{x^{1+\alpha_p}}\mathrm{e}^{-\lambda_p x}1_{x>0}\,\mathrm{d}x +
\frac{c_n}{(-x)^{1+\alpha_n}}\mathrm{e}^{\lambda_n x}1_{x<0}\,\mathrm{d}x,
\]
where $\alpha_p, \alpha_n <2$, $\lambda_p, \lambda_n, c_p, c_n >0$.
They satisfy (\ref{levy_cond}) for any $\alpha< \min(\lambda_p,
\lambda_n)$ and have infinite variation when $\alpha_p, \alpha_n
\geq1$. Another example are Meixner processes, used in financial
modeling (cf. Kyprianou et al. \cite{KSW}, Chapter I), with L\'{e}vy
measure given by
\[
\nu(\mathrm{d}x) = \delta\frac{\exp (\beta x/ \eta )}{x \sinh
( \uppi x/\eta )}\,\mathrm{d}x,
\]
where $\delta, \eta>0$, $|\beta| < \uppi$. They satisfy (\ref
{levy_cond}) for $\alpha<  (\uppi- |\beta| )/\eta$.

\section*{Acknowledgements}
The authors would like to thank the anonymous referee and the Associate
Editor for many valuable comments and suggestions, which considerably
improved the presentation of the results. The research of the first and
the second author was supported by the National Science Center in
Poland under decision no.
DEC-2011/01/B/ST1/05089.




\printhistory

\end{document}